\documentclass[12pt,reqno,oneside]{amsart}

\usepackage[english]{babel}
\usepackage{amsmath}
\usepackage{amssymb, amsthm}
\usepackage{enumerate}
\usepackage{comment}
\usepackage[all,cmtip]{xy}
\usepackage[bottom=3.5cm]{geometry}
\usepackage{graphicx}
\usepackage{array}
\usepackage{url}
\usepackage{hyperref}
\usepackage{setspace}
\usepackage[dvipsnames]{xcolor}
\usepackage{tikz}
\usepackage{wrapfig,lipsum}
\usepackage{todonotes}
\usepackage{etex}
\usepackage{flafter}
\usetikzlibrary{shapes}
\hypersetup{
  linktoc=page,
}
\author{Markus Hausmann}
\date{}
\newtheorem{satz}{Satz}[section]
\newtheorem{Theorem}[satz]{Theorem}
\newtheorem{Cor}[satz]{Corollary}
\newtheorem{Lemma}[satz]{Lemma}
\newtheorem{Prop}[satz]{Proposition}

\theoremstyle{definition}
\newtheorem{Def}[satz]{Definition}

\newtheorem{Remark}[satz]{Remark}
\newtheorem{Example}[satz]{Example}
\newtheorem*{Def*}{Definition}
\newtheorem*{Example*}{Example}
\newtheorem*{Theorem*}{Theorem}
\newtheorem*{Prop*}{Proposition}
\newtheorem*{Remark*}{Remark}

\newcommand{\R}{\mathbb{R}}
\newcommand{\Q}{\mathbb{Q}}
\newcommand{\Z}{\mathbb{Z}}
\newcommand{\cC}{\mathcal{C}}

\newcommand{\upi}{\underline{\pi}}
\newcommand{\N}{\mathbb{N}}
\newcommand{\F}{\mathcal{F}}
\newcommand{\U}{\mathcal{U}}
\newcommand{\walph}{\widetilde{\alpha}}
\newcommand{\lr}{\overline{r}}
\newcommand{\xr}{\xrightarrow}

\newcommand{\dia}{\diamond}
\newcommand{\wL}{\widetilde{L}}

\newcommand{\mC}{\mathcal{C}}

\newcommand{\rrho}{\overline{\rho}}

\DeclareMathOperator*{\colim}{colim}
\DeclareMathOperator*{\hocolim}{hocolim}
\DeclareMathOperator*{\shift}{sh}

\DeclareMathOperator*{\tr}{tr}

\DeclareMathOperator*{\Aut}{Aut}

\DeclareMathOperator*{\im}{im}

\DeclareMathOperator*{\Epi}{Epi}
\DeclareMathOperator*{\Out}{Out}
\DeclareMathOperator*{\Ch}{Ch}
\DeclareMathOperator*{\Sing}{Sing}

\makeatletter
\newcommand*{\defeq}{\mathrel{\rlap{%
                     \raisebox{0.3ex}{$\m@th\cdot$}}%
                     \raisebox{-0.3ex}{$\m@th\cdot$}}%
                     =}
\makeatother

\tikzset{
    partial ellipse/.style args={#1:#2:#3}{
        insert path={+ (#1:#3) arc (#1:#2:#3)}
    }
}

\numberwithin{equation}{section}

\title{Symmetric products and subgroup lattices}

\begin{document}

\begin{abstract}Let $G$ be a finite group. We show that the rational equivariant homotopy groups of symmetric products of the $G$-equivariant sphere spectrum are naturally isomorphic to the rational homology groups of certain subcomplexes of the subgroup lattice of~$G$.
\end{abstract}
\maketitle

\section{Introduction}
Let $Sp^n(X)=X^{\times n}/\Sigma_n$ denote the $n$-th symmetric product of a space $X$, and $Sp^n=\{Sp^n(S^k)\}$ the spectrum consisting of the $n$-th symmetric products of spheres. Inserting a basepoint in the last component defines maps $Sp^n\to Sp^{n+1}$, giving rise to the symmetric product filtration
\[ \mathbb{S}=Sp^1\to Sp^2 \to \hdots \to Sp^{\infty}\simeq H\Z \]
which interpolates between the sphere spectrum and the Eilenberg-MacLane spectrum for the integers. This filtration has interesting properties: For example, it induces the filtration by length of admissible sequences on the Steenrod algebra (Nakaoka \cite{Nak58}), it is related to partition complexes and the Goodwillie tower of the identity (Arone-Dwyer \cite{AD01}) and it is the object of study in the Whitehead conjecture (solved by Kuhn \cite{Kuh82}). However, all these properties are purely about torsion: After tensoring with~$\Q$, the symmetric product filtration becomes constant.

Now let $G$ be a finite group. The symmetric products $Sp^n(S^V)$ of $G$-representation spheres give rise to a genuine $G$-spectrum, which we denote by $Sp^n_G$. Again one obtains a filtration 
\[ \mathbb{S}_G=Sp_G^1\to Sp_G^2 \to \hdots \to Sp_G^{\infty}\simeq H\underline{\Z}, \]
this time converging to an Eilenberg-MacLane spectrum for the constant Mackey functor~$\underline{\Z}$ (dos Santos \cite{dS03}). In contrast to the non-equivariant situation, the Hurewicz map $\mathbb{S}_G\to Sp^{\infty}_G$ is no longer a rational equivalence, as one can see on the level of $\pi^G_0$: The group $\pi_0^G(\mathbb{S}_G)\otimes \Q$ is the rationalized Burnside ring, a vector space with basis the isomorphism classes of transitive $G$-sets, while $\pi_0^G(Sp^{\infty}_G)\otimes \Q$ is isomorphic to $\Q$. It is still the case, however, that all higher rational homotopy groups both of $\mathbb{S}_G$ and $Sp^{\infty}_G$ vanish.

The content of this paper is to show that the intermediate rational homotopy groups $\pi_k^G(Sp_G^n)\otimes \Q$ are closely related to the topology of the subgroup lattice of $G$, and in general not at all concentrated in degree~$0$. In fact, the only $G$ for which all $\pi_k^G(Sp_G^n)\otimes \Q$ with $k>0$ vanish are the cyclic $p$-groups. If one lets $G$ and~$n$ vary, the vector spaces $\pi_k^G(Sp_G^n)\otimes \Q$ can be non-trivial for arbitrarily large $k$, and for every $k$ these can be of arbitrarily large finite dimension.

The precise statement is as follows: We recall that the (nerve of the) subgroup lattice $L(G)$ is the simplicial set with $k$-simplices the chains of inclusions $H_0\leq \hdots \leq H_k$ of subgroups of $G$. We define a filtration
\[ \emptyset=L(G)_0\subseteq L(G)_1\subseteq L(G)_2 \subseteq \hdots \subseteq L(G)_{\infty}=L(G) \]
on this subgroup lattice by declaring a simplex $H_0\leq \hdots \leq H_k$ to lie in $L(G)_n$ if and only if the index $[H_k:H_0]$ is at most $n$. The lattice $L(G)$ carries an action by $G$ via conjugation, which preserves the subcomplexes~$L(G)_n$. Hence, the rational homology $H_*(L(G)_n,\Q)$ becomes a graded $G$-module and we can consider its co-invariants $(H_*(L(G)_n,\Q))_G$. Then the main result of this paper is the following:
\begin{Theorem} \label{thm:intro} For all finite groups $G$ and $n\in \N\cup \{\infty\}$ there are isomorphisms
\[ \pi_*^G(Sp_G^n)\otimes \Q \cong (H_*(L(G)_n,\Q))_G. \]
\end{Theorem}
In other words, the process of adding the $n$-th coordinate in the symmetric products of $G$-spheres can be modeled rationally by adding all chains of total index $n$ in the subgroup lattice of $G$. In particular, the rational homotopy type of $Sp_G^n$ only changes when $n$ is a divisor of the order of $G$. We observe that Theorem \ref{thm:intro} matches the previously known values $\pi_*^G(\mathbb{S}_G)\otimes \Q$ and $\pi_*^G(Sp^{\infty}_G)\otimes \Q$, since $L(G)_1$ is the discrete set of subgroups of $G$ and $L(G)_{\infty}=L(G)$ is contractible since the lattice has a minimal and a maximal element. As mentioned above, the subcomplexes $L(G)_n$ can have arbitrarily high non-trivial rational homology, though of course for every fixed $G$ the homology is bounded since $L(G)$ is a finite complex.
For groups of small order, Theorem \ref{thm:intro} makes it an easy exercise to determine each $\pi_*^G(Sp^n_G)\otimes \Q$ concretely, and we work this out in several examples. For general $G$, however, the rational homology of the $L(G)_n$ is difficult to compute and Theorem \ref{thm:intro} can be seen as a demonstration that the rationalized $Sp_G^n$ are quite complex, in contrast to their non-equivariant version.

\begin{Remark} Schwede \cite{Sch17} has given a computation of the $0$-th homotopy groups $\pi_0^G(Sp^n_G)$, even integrally. His description does not make use of the filtration of subgroup lattices that we introduce in this paper. Instead, he showed that each $\pi_0^G(Sp^n_G)$ is a quotient of $\pi_0^G(\mathbb{S})$, the Burnside ring of $G$, and that the kernel is generated by a single element when viewed as a global functor (see the discussion below). It is not hard to see that -- after tensoring with $\Q$ -- his description agrees with ours. One can think of Theorem \ref{thm:intro} as a higher homotopical generalization of the rational version of Schwede's result.
\end{Remark}

To describe the full functoriality in $G$, we work in the global equivariant context, as introduced by Schwede in \cite{Sch18}. This means that we think of the collection of all $Sp^n_G$, for varying $G$ and fixed $n$, as one global object. One consequence is that the homotopy groups form a global functor: They carry restriction maps along arbitrary group homomorphisms and transfer maps for subgroup inclusions. The strongest form of our main result (Theorem \ref{theo:strong}) also describes the global functor structure on the side of subgroup lattices, and in addition explains how to reconstruct the full global homotopy type of the rationalization $Sp^n_{\Q}$ in terms of subgroup lattice data, not only its homotopy groups. This makes a difference, because - unlike for fixed $G$ - a rational global spectrum is not determined by its homotopy groups. In fact, this stronger version allows us to deduce that the $Sp^n_{\Q}$ do not decompose as products of Eilenberg-MacLane spectra globally (Proposition \ref{prop:formal}), though they decompose over every fixed $G$. To obtain these results we make use of an explicit construction of an equivalence between the homotopy category of rational global spectra and the derived category of rational global functors, due to Wimmer \cite{Wim17}.

\begin{Remark} One ingredient in the proof of Theorem \ref{thm:intro} is an equivariant version of the theorem by Arone and Dwyer \cite[Theorem 1.11]{AD01} which relates $Sp^n/Sp^{n-1}$ to the partition complex $\Pi_n$ of the set $\{1,\hdots,n\}$. This is Proposition \ref{prop:ad} and might be of independent interest. 
\end{Remark}

The \textbf{organization} of the paper is as follows: In Section \ref{sec:ratglobal} we recall basics of (rational) global homotopy theory, explain how Theorem \ref{thm:intro} can be stated in terms of geometric instead of categorical fixed points and formulate the stronger version via chain complexes of global functors.
In Sections \ref{sec:geomidea} - \ref{sec:construction} we construct compatible maps from the suspension spectrum of the $L(G)_n$ into the geometric fixed points of $Sp^n$. The basic geometric idea is not complicated, but it takes some effort to obtain an honest map with all the properties and compatibilities we need. We then show (Section \ref{sec:split}) that the induced map on subquotients is rationally (split) injective on homology. In Section \ref{sec:ad} we use other methods to compute the rational homotopy type of the geometric fixed points of $Sp^n/Sp^{n-1}$ and show that it agrees with the rational homotopy type of the lattice quotient, finishing the proof of the main results.
In Section \ref{sec:examples} we give some examples for small~$G$. Finally, Section \ref{sec:globalprop} contains a proof that the $Sp^n_{\Q}$ do not split as products of global Eilenberg-MacLane spectra, and that the $\pi_*^G(Sp^n)\otimes \Q$ are only concentrated in degree $0$ when $G$ is a cyclic $p$-group.

\textbf{Acknowledgements}: I would like to thank my adviser Stefan Schwede for various helpful discussions and Christian Wimmer for answering my questions on rational global homotopy theory.

This research was supported by the GRK 1150 `Homotopy and Cohomology' and the DFG Priority Program SPP 1786 `Homotopy Theory and Algebraic Geometry'. Final revisions were made in Copenhagen under the support of the Danish National Research Foundation through the DNRF92 `Centre for Symmetry and Deformation'.

\section{Rational global homotopy} \label{sec:ratglobal}

In this section we recall some definitions and structural results of (rational) stable global homotopy theory based on orthogonal spectra, as developed by Schwede in his book project \cite{Sch18} and the paper \cite{Sch17}. Here and throughout the paper, we work with the version formed with respect to the global family of \emph{finite} groups, not the full one for all compact Lie groups. In particular, we explain how one can reduce our main result to a statement about geometric fixed points. For the reader only interested in the version for fixed $G$, we recall in Remarks \ref{rem:ggeom} and \ref{rem:gmorita} how this reduction works in the more familiar category of $G$-spectra.

\subsection{Orthogonal spectra} An orthogonal spectrum is a collection of based spaces $X(V)$ for every finite dimensional real inner product space $V$, together with associative structure maps of the form $X(V)\wedge S^{W-\varphi(V)}\to X(W)$, varying continuously in linear isometric embeddings $\varphi:V\hookrightarrow W$ (see Mandell-May-Schwede-Shipley \cite[Example I.4.4]{MMSS01}, Schwede \cite[Section 2]{Sch17} or Hill-Hopkins-Ravenel \cite[Section A.2.4]{HHR16} on how to make this precise). In particular, each $X(V)$ carries a based $O(V)$-action. If an inner product space $V$ comes with an action of a finite group $G$, the evaluation $X(V)$ also becomes a based $G$-space through functoriality. If $\varphi:V\hookrightarrow W$ is an equivariant embedding of such representations, the structure map $X(V)\wedge S^{W-\varphi(V)}\to X(W)$ becomes $G$-equivariant. This way we think of an orthogonal spectrum as a collection of $G$-spaces for all finite groups $G$ which are related by structure maps, i.e., a global spectrum. In particular, every orthogonal spectrum $X$ gives rise to a genuine $G$-orthogonal spectrum (in the sense of Mandell-May \cite{MM02}) $X_G$ by only remembering the evaluations at $G$-representations.

\begin{Example}[Symmetric products] In this paper we study the orthogonal spectra~$Sp^n$. Their evaluation on $V$ is given by $Sp^n(S^V)=(S^V)^{\times n}/\Sigma_n$. For a linear isometric embedding $\varphi:V\hookrightarrow W$, the associated structure map
\[ (S^V)^{\times n}/\Sigma_n\wedge S^{W-\varphi(V)}\to (S^W)^{\times n}/\Sigma_n \]
sends a class $[(v_1,\hdots,v_n)]\wedge x$ to $[(\varphi(v_1)\wedge x,\hdots, \varphi(v_n)\wedge x)]$. The insertion of a basepoint gives inclusions $i_{n-1}^n:Sp^{n-1}\to Sp^n$ with colimit spectrum $Sp^{\infty}$. More generally, for $m\leq n$ we write $i_m^n:Sp^m\to Sp^n$ for the map obtained by inserting $(n-m)$ basepoints.
\end{Example}
For every orthogonal spectrum $X$ there are two notions of equivariant homotopy groups, the categorical fixed point homotopy groups $\pi_*^G(X)$ (which Theorem \ref{thm:intro} is about) and the geometric fixed point homotopy groups $\Phi_*^G(X)$. Rationally, the two determine each other. In the following sections we quickly recall their definition and how to pass from one to the other.

\subsection{Categorical homotopy groups} \label{sec:cathom}
Let $\rho_G$ denote the regular representation of a finite group $G$.
\begin{Def} The categorical homotopy groups of an orthogonal spectrum $X$ are defined as
\[ \pi_k^G(X)=\colim_{n\in \N}[S^{k+n\cdot \rho_G},X(n\cdot \rho_G)]^G, \]
for all finite groups $G$ and all $k\in \Z$, where $[-,-]^G$ denotes the set of based homotopy classes of based $G$-maps. The colimit is formed with respect to the structure maps $X(n\cdot \rho_G)\wedge S^{\rho_G}\to X((n+1)\cdot \rho_G)$.
\end{Def}

These groups have the following functoriality in $G$ (cf. Schwede \cite[Constructions 3.2.22 and 3.1.15]{Sch18}):
\begin{itemize}
	\item For every subgroup inclusion $H\leq G$ there is a \emph{transfer} map 
	\[ {\tr}_H^G:\pi_*^H(X)\to \pi_*^G(X). \]
	\item For every group homomorphism $\psi:G\to K$ there is a \emph{restriction} map \[ \psi^*:\pi_*^K(X)\to \pi_*^G(X). \]
	If $K=G$ and $\psi$ is an inner automorphism, then $\psi^*$ is the identity of $\pi_*^G(X)$.
\end{itemize}
Transfer maps are covariantly functorial in the group $G$, and restriction maps are contravariantly functorial. Restrictions along surjective group homomorphisms commute with transfers - in the appropriate sense - while restrictions along injective group homomorphisms are related to transfers by a double coset formula. The precise formulations can be found in \cite[Theorem 4.2.6 and the paragraph following it]{Sch18}. Altogether, the collection $\upi_*(X)=\{\pi_*^G(X)\}_{G\text{ finite}}$ forms a so-called graded \emph{global functor}. In \cite{Sch18}, Schwede defines a morphism of orthogonal spectra to be a \emph{global equivalence} if it induces an isomorphism of graded global functors on homotopy groups. The \emph{global homotopy category} is the localization of the category of orthogonal spectra at the class of global equivalences.

Likewise, a \emph{rational global equivalence} is a morphism that induces an isomorphism on $\upi_*(-)\otimes \Q$. The localization of orthogonal spectra at rational global equivalences forms the \emph{rational global homotopy category}. As in the non-equivariant case, the passage from the global homotopy category to the rational global homotopy category is a left Bousfield localization. It has a fully faithful right adjoint with essential image those orthogonal spectra whose categorical homotopy groups form $\Q$-vector spaces. Given an orthogonal spectrum~$X$, its rationalization $X_{\Q}$ can be constructed as the homotopy colimit of the sequence $X\xr{\cdot 2}X\xr{\cdot 3}X\xr{\cdot 4} \hdots$, or as the smash product with the rational global sphere $\mathbb{S}_{\Q}$. Furthermore, the counit $X\to X_{\Q}$ induces an isomorphism $\upi_*(X)\otimes \Q\cong \upi_*(X_{\Q})$.

Now we let $L(G)_n$ denote the subcomplex of the subgroup lattice $L(G)$ consisting of all chains of total index at most $n$, as explained in the introduction. Since we aim to show that $(H_*(L(G)_n,\Q))_G$ is isomorphic to $\pi^G_*(Sp^n)\otimes \Q$ for every finite group $G$, the assignment
\[ G\mapsto (H_*(L(G)_n,\Q))_G \]
must have the functoriality of a global functor. In fact, this structure already exists on rational chains, i.e., the assignment
\[ G\mapsto (\mathcal{C}_*(L(G)_n))_G \]
extends to a chain complex of rational global functors, where $\mathcal{C}_*$ denotes the normalized rational complex of a simplicial set. The transfer and restriction maps work out as follows:

\textbf{Transfers}: The transfer along a subgroup inclusion $H\leq G$ sends a chain of subgroups of $H$ to $[G:H]$ times the same chain thought of as subgroups of~$G$.

\textbf{Restriction maps}: The restriction along a group homomorphism $\psi:G\to K$ takes the class of a chain~$H_0\leq \hdots \leq H_m$ to
\begin{equation} \label{eq:doublecoset} \sum_{[k]\in G_{\psi}\backslash K/H_0} \left(\frac{[G:\psi^{-1}(kH_0k^{-1})]}{[K:H_0]}\cdot [\psi^{-1}(kH_0k^{-1})\leq \hdots \leq \psi^{-1}(kH_mk^{-1})]\right). \end{equation}
The sum is taken over a set of coset representatives of the $(G\times H_0^{op})$-action on $K$ given by $(g,h)\cdot k=\psi(g)kh$. We note that if $\psi$ is surjective, there is only one summand and the formula simplifies to
\[ \psi^*([H_0\leq \hdots \leq H_m])=[\psi^{-1}(H_0)\leq \hdots \leq \psi^{-1}(H_m)]. \]

\begin{Remark} \label{rem:fraction} To explain the factor $[G:H]$ in the definition of the transfer and the fractions $\frac{[G:\psi^{-1}(kH_0k^{-1})]}{[K:H_0]}$ in that of the restrictions, we consider the case $n=1$, where $L(G)_1$ is the discrete set of subgroups of $G$ and $\pi_0^G(Sp^1)\otimes \Q$ is the rationalized Burnside ring of $G$, a vector space with basis the isomorphism classes of transitive $G$-sets $G/H$. It is clear that $(H_0(L(G)_1,\Q))_G$ and $\pi_0^G(Sp^1)\otimes \Q$ are abstractly isomorphic, and one isomorphism would be given by simply mapping $H\in L(G)_1$ to $[G/H]\in \pi_0^G(Sp^1)\otimes \Q$. However, this isomorphism cannot be compatible with any choice of isomorphism 
\[ (H_0(L(G)_{\infty},\Q))_G\cong \pi_0^G(Sp^{\infty})\otimes \Q, \]
which are both isomorphic to $\Q$. Any two subgroups of $G$ represent the same element in $H_0(L(G)_{\infty},\Q)$, since the lattice is connected. On the other hand, the augmentation $\pi_0^G(Sp^1)\to \pi_0^G(Sp^{\infty})$ sends a finite $G$-set to its number of elements, so the orbits $[G/H]$ generally have different images. This can be corrected by sending $H\in L(G)_1$ to $\frac{1}{[G:H]}\cdot [G/H]$ - which has augmentation $1$ - instead. So, in the case $n=1$, the formulas above express transfers and restrictions in the rationalized Burnside ring global functor written in terms of the basis $\{\frac{1}{[G:H]}\cdot [G/H]\}$ instead of the usual basis $\{[G/H]\}$, which leads to the appearance of the factors and fractions.
\end{Remark}

We claim that these transfers and restrictions turn $(\Q[L(-)_n])_-$ into a simplicial global functor. The proof that each simplicial degree is a global functor is very similar to the Burnside ring global functor (the case $n=1$, see Remark \ref{rem:fraction} above), since the main role in the restrictions \eqref{eq:doublecoset} is played by the smallest subgroup $H_0$ and the higher $H_i$ are carried along. The proof that the global structure maps commute with the simplicial operators is straightforward, except possibly for the face $d_0$ and a restriction~$\psi$. There we have to show that the two sums
\begin{equation} \label{eq:form1} \sum_{[k]\in G_{\psi}\backslash K/H_0} \left(\frac{[G:\psi^{-1}(kH_0k^{-1})]}{[K:H_0]}\cdot [\psi^{-1}(kH_1k^{-1})\leq \hdots \leq \psi^{-1}(kH_mk^{-1})]\right) \end{equation}
and
\begin{equation} \label{eq:form2} \sum_{[\widetilde{k}]\in G_{\psi}\backslash K/H_1} \left(\frac{[G:\psi^{-1}(\widetilde{k}H_1\widetilde{k}^{-1})]}{[K:H_1]}\cdot [\psi^{-1}(\widetilde{k}H_1\widetilde{k}^{-1})\leq \hdots \leq \psi^{-1}(\widetilde{k}H_m\widetilde{k}^{-1})]\right) \end{equation}
are the same. The class $[\psi^{-1}(kH_1k^{-1})\leq \hdots \leq \psi^{-1}(kH_mk^{-1})]$ only depends on the $(G\times H_1^{op})$-orbit of $k$, and so we can rewrite \eqref{eq:form1} as
\[ \sum_{\substack{[\widetilde{k}]\in \\ G_{\psi}\backslash K/H_1}}\left(\left(\sum_{\substack{[k]\in G_{\psi}\backslash K/H_0\\ k\in G_{\psi}\widetilde{k}H_1}} \frac{[G:\psi^{-1}(kH_0k^{-1})]}{[K:H_0]}\right)\cdot [\psi^{-1}(\widetilde{k}H_1\widetilde{k}^{-1})\leq \hdots \leq \psi^{-1}(\widetilde{k}H_m\widetilde{k}^{-1})]\right).
\]
Hence, it suffices to show that
\[  \sum_{\substack{[k]\in G_{\psi}\backslash K/H_0\\ k\in G_{\psi}\widetilde{k}H_1}} \frac{[G:\psi^{-1}(kH_0k^{-1})]}{[K:H_0]} = \frac{[G:\psi^{-1}(\widetilde{k}H_1\widetilde{k}^{-1})]}{[K:H_1]} \] 
for every $\widetilde{k}\in K$. This equality can be deduced from counting the number of elements in the $(G\times H_1^{op})$-orbit of $\widetilde{k}\in K$ in two different ways: The $(G\times H_1^{op})$-isotropy of $\widetilde{k}$ is given by the graph subgroup
\[ \{(g,\widetilde{k}^{-1} \psi(g)^{-1} \widetilde{k})\ |\ g\in \psi^{-1}(\widetilde{k} H_1 \widetilde{k}^{-1}) \}, \]
so the order of this orbit is \[ \frac{|G|\times |H_1|}{|\psi^{-1}(\widetilde{k} H_1 \widetilde{k}^{-1})|}.\]
On the other hand, decomposing the orbit $G_{\psi}\widetilde{k}H_1$ into $(G\times H_0^{op})$-orbits yields the sum 
\[ \sum_{\substack{[k]\in G_{\psi}\backslash K/H_0\\ k\in G_{\psi}\widetilde{k}H_1}} \frac{|G|\times |H_0|}{|\psi^{-1}(k H_1 k^{-1})|}, \]
and so dividing by $|K|$ gives the desired result.

Hence, the normalized chains associated to $(\Q[L(-)_n])_-$ become a complex of global functors, which we denote by~$\cC L_n$. Then a more functorial version of Theorem \ref{thm:intro} is the following:
\begin{Theorem} \label{thm:mackey} There is an isomorphism of graded global functors
\[ \upi_* (Sp^n)\otimes \Q\cong H_*(\cC L_n). \]
\end{Theorem}

\subsection{Geometric homotopy groups} \label{sec:geom} Next we explain how to reduce Theorem \ref{thm:mackey} to a statement about geometric fixed point homotopy groups, which in the case of symmetric products of spheres are more directly accessible.
\begin{Def}[Geometric fixed points] The geometric fixed point homotopy groups of an orthogonal spectrum $X$ are defined as
 \[ \Phi_k^G(X)=\colim_{n\in \N} [S^{k+n},X(n\cdot \rho_G)^G] \]
for every finite group $G$ and $k\in \Z$.
\end{Def}
The collection of geometric fixed point homotopy groups does not form a global functor, but it carries restrictions along surjective group homomorphisms (\cite[Construction 3.3.4]{Sch18}). Again, inner conjugations act trivially. In other words, if we denote the category of finite groups and conjugacy classes of surjective group homomorphisms by $\Out$, the collection $\Phi_*(X)$ forms a functor $\Out^{op}\to Ab$, an $\Out^{op}$-module (\cite[Proposition 4.1.23]{Sch18}). Moreover, a morphism of orthogonal spectra is a global equivalence if and only if it induces an isomorphism on all geometric fixed point homotopy groups.

There is a natural comparison map $\gamma_X:\pi_k^G(X)\to \Phi_k^G(X)$ given by taking a $G$-map $S^{k+n\cdot \rho_G}\to X(n\cdot \rho_G)$ to the induced map on fixed points $S^{k+n}\cong (S^{k+n\cdot \rho_G})^G\to X(n\cdot \rho_G)^G$. The map $\gamma_X$ has the following properties:
\begin{enumerate}[(i)]
\item It commutes with restrictions along surjective group homomorphisms.
\item It takes all elements of the form $\tr_H^G(x)$ for $H$ a proper subgroup of $G$ to $0$.
\end{enumerate}
Given a global functor $F$ and a finite group $G$ we let $\tau(F)(G)$ denote the quotient of $F(G)$ by all transfers from proper subgroups. Then the assignment $G\mapsto \tau(F)(G)$ no longer forms a global functor, but it inherits restrictions along surjective group homomorphisms, since these commute with transfers. In these terms, the two properties above mean that $\gamma_X$ factors through a map of $\Out^{op}$-modules $\widetilde{\gamma}_X:\tau(\upi_*(X))\to \Phi_*(X)$. It turns out that for rational global functors the construction $\tau$ can be reversed and that in this case $\widetilde{\gamma}_X$ is an isomorphism:
\begin{Prop}[Schwede {\cite[Theorem 4.5.35 and Corollary 4.5.37]{Sch18}}] \label{prop:equiv} The functor
\begin{align*} \tau:\Q-\textnormal{global functors} & \to  \Q[{\Out}^{op}]-\textnormal{mod}
\end{align*}
is an equivalence of categories. Moreover, for every orthogonal spectrum $X$ the map
\[\widetilde{\gamma}_X:\tau(\upi_*(X)\otimes \Q)\to \Phi_*(X)\otimes \Q \]
is an isomorphism.
\end{Prop}
\begin{Example} In case of the global sphere spectrum $\mathbb{S}$, the geometric fixed points $\Phi^G(\mathbb{S})$ for any finite group~$G$ are given by the non-equivariant sphere spectrum. The map $\gamma_{\mathbb{S}}:\pi^G_0(\mathbb{S})\otimes \Q\to \Phi_0^G(\mathbb{S})\otimes \Q\cong \Q$ sends the class of a finite $G$-set $X$ (as an element in the Burnside ring) to the cardinality of its fixed points $X^G$, or in other words the number of trivial orbits in a decomposition $X\cong \bigsqcup G/H_i$. This map is a surjection, with kernel generated by all transitive $G$-sets $G/H$ for proper subgroups $H$. These also span the subspace generated by transfers, since each $G/H$ is the transfer of the trivial $H$-set $H/H$. Using that all higher rational homotopy groups of $\mathbb{S}$ vanish, we see that $\widetilde{\gamma}_{\mathbb{S}}$ is indeed an isomorphism.
\end{Example}
\begin{Remark} \label{rem:ggeom} When working over a fixed finite group $G$, the analogous statement is the equivalence of categories between rational $G$-Mackey functors and products of $\Q[W_G(H)]$~-modules, where $H$ ranges through a set of conjugacy class representatives and $W_G(H)=N_G(H)/H$ denotes the Weyl group of $H$, cf. \cite[Appendix A]{GM95} or \cite[Theorem 3.4.22]{Sch18}.
\end{Remark}
This implies that instead of showing that $\upi_*(Sp^n)\otimes \Q$ is isomorphic to $H_*(\cC L_n)$ as a global functor (Theorem \ref{thm:mackey}), we can equivalently show:
\begin{Theorem} \label{thm:geom} There is an isomorphism of graded $\Out^{op}$-modules
\[ \Phi_*(Sp^n)\otimes \Q\cong \tau(H_*(\cC L_n)). \]
\end{Theorem}
The graded $\Out^{op}$-module $\tau(H_*(\cC L_n))$ can be described as follows. Since $\tau$ is an equivalence, it is in particular exact. So there is an isomorphism
\[ \tau(H_*(\cC L_n))\cong H_*(\tau(\cC L_n)). \]
By the description given in Section \ref{sec:cathom}, the transfers of $\cC L_n(G)$ are generated by all chains of subgroup inclusions that end in a proper subgroup of $G$. This process can be carried out on the space level: Let $\widetilde{L}(G)$ denote the quotient of $L(G)$ by the subspace of chains that do not end in $G$, and similarly define $\widetilde{L}(G)_n$. Then we find that there are isomorphisms
\[ \tau(\cC L_n)(G)\cong (\mathcal{C}_*(\widetilde{L}(G)_n))_G\cong \mathcal{C}_*(\widetilde{L}(G)_n/G), \]
for every finite group $G$. Under this isomorphism, the restriction along a surjection $\psi:G\twoheadrightarrow K$ sends a chain $[H_0\leq \hdots \leq H_m]$ to $[\psi^{-1}(H_0)\leq \hdots \leq \psi^{-1}(H_m)]$. In other words, the $\Out^{op}$-complex $\tau(\cC L_n)$ arises by applying rational chains to the simplicial set valued $\Out^{op}$-functor $G\mapsto \widetilde{L}(G)_n/G$. For this reason we from now on write
$\cC \widetilde{L}_n$ instead of $\tau(\cC L_n)$.

\subsection{Global chain complexes} \label{sec:globalchain}
In fact we show something stronger than Theorem \ref{thm:geom}, for which we recall some more rational global homotopy theory.
As an application of Morita theory for stable model categories (see \cite{SS03}), the model category of orthogonal spectra with rational global equivalences is Quillen equivalent to the derived category of rational global functors (\cite[Theorem 4.5.29]{Sch18}). The latter in turn -- using Proposition \ref{prop:equiv} above -- is equivalent to the derived category of rational $\Out^{op}$-modules. In his PhD thesis \cite{Wim17}, Wimmer constructs an explicit equivalence $T$ between the rational global homotopy category and this derived category. We now recall his construction.

Let $\Epi$ denote the category of finite groups and surjective group homomorphisms. Then the equivalence is constructed as a composite
\[ \text{orthogonal spectra} \xr{\Phi} (\text{orthogonal spectra})^{\Epi^{op}} \xr{c^{\Epi^{op}}} {\Ch}_{\Q}^{\Epi^{op}} \xr{q_{!}} {\Ch}_{\Q}^{\Out^{op}}, \]
where $\Ch_{\Q}$ is the category of rational chain complexes and $ (\text{orthogonal spectra})^{\Epi^{op}}, \newline {\Ch}_{\Q}^{\Epi^{op}}$ and ${\Ch}_{\Q}^{\Out^{op}}$ denote the respective functor categories.

We go through each of these functors individually:
The first, $\Phi$, sends an orthogonal spectrum $X$ to the collection of geometric fixed point spectra $\Phi^G(X)$ defined via
\[ \Phi^G(X)(V)=X(V\otimes \rho_G)^G.\]
Given a surjection $\psi:G\twoheadrightarrow K$, the induced morphism 
\[ \psi^*:\Phi^K(X)\to \Phi^G(X)\] is defined in level $V$ by
\[ X(V\otimes \rho_K)^K= X(V\otimes \psi^*(\rho_K))^G\xr{X(V\otimes i_{\psi})^G} X(V\otimes \rho_G)^G. \] 
Here, $\psi^*(\rho_K)$ denotes the restriction of $\rho_K$ to a $G$-representation along $\psi$, and $i_{\psi}:\psi^*(\rho_K)\hookrightarrow \rho_G$ is the $G$-equivariant linear isometry which sends a basis element $e_k$ to 
\[ \sqrt{\frac{|K|}{|G|}}\sum_{g\in \psi^{-1}(k)} e_g.\]
To avoid confusion, one should no longer think of $\Phi^G(X)$ as a global spectrum. Only the non-equivariant homotopy type is important. Almost by definition, there is a natural isomorphism between $\Phi_*^G(X)$ and the non-equivariant homotopy groups $\pi_*(\Phi^G(X))$. Moreover, the restriction maps on $\Phi(X)$ induce the restriction maps of $\Phi_*(X)$ under this isomorphism.

The second functor $c^{\Epi^{op}}$ is given by postcomposition with a functor 
\[ c:\text{orthogonal spectra}\to {\Ch}_{\Q}, \]
of which we need the following two properties:
\begin{itemize}
	\item There is a natural isomorphism between $\pi_*(X)\otimes \Q$ and $H_*(c(X))$. Hence, $c$ takes rational equivalences to quasi-isomorphisms.
	\item For all based spaces $A$ there is a natural quasi-isomorphism $\cC_*(A)\simeq c(\Sigma^{\infty}A)$. In other words, for suspension spectra of spaces the associated rational chain complex $c(\Sigma^{\infty}A)$ is equivalent to the usual rational singular chains.
\end{itemize}
We quickly give the construction of such a functor. Wimmer uses a more refined version that also preserves multiplicative structures, but the one we describe now is enough for our purposes. Let $X$ be an orthogonal spectrum. Define $c(X)$ as the sequential colimit
\[ \colim_{n\in \mathbb{N}} \cC_*(X(\mathbb{R}^n))[-n], \]
where again $\cC_*(-)$ stands for the rational chains on a based space and $[-n]$ denotes the $(-n)$-fold shift of a chain complex, i.e., $(C[-n])_k=C_{k+n}$ and $d^{C[-n]}_k=d^C_{k+n}$. The connecting maps are given by the composite
\[  \cC_*(X(\mathbb{R}^n))[-n]\xrightarrow{\nabla (- \otimes \iota)} \cC_*(X(\mathbb{R}^n)\wedge S^1)[-(n+1)] \xrightarrow{(\sigma_{\mathbb{R}^n}^{\mathbb{R}})_*}\cC_*(X(\mathbb{R}^{n+1}))[-(n+1)]. \]
The first map is the Eilenberg-Zilber shuffle product with a fixed integral $1$-cycle $\iota$ of $S^1$ that represents a generator in $H_1(S^1,\mathbb{Z})$, and the second map is induced by the structure map of $X$. That the homology of $c(X)$ computes the rational homology (and hence the rational homotopy) of $X$ follows from the fact that $X$ is equivalent to the homotopy colimit over the spectra $\Sigma^{-n} (\Sigma^{\infty}X(\mathbb{R}^n))$ (which is standard for sequential spectra and also holds for orthogonal spectra, since the forgetful functor to sequential spectra preserves stable equivalences). In case $X$ is a suspension spectrum $\Sigma^{\infty}A$, all structure maps $\sigma_{\mathbb{R}^n}^{\mathbb{R}}$ are homeomorphisms. Since the Eilenberg-Zilber shuffle product with $\iota$ is a quasi-isomorphism, it follows that all connecting maps in the colimit system are quasi-isomorphisms and so the canonical natural map $\cC_*(A)=\cC_*((\Sigma^{\infty}A)(0))\to c(\Sigma^{\infty}A)$ is also one, which proves the second desired property.

Finally, $q:\Epi \to \Out$ is the projection and we write $q_{!}:{\Ch}_{\Q}^{\Epi^{op}} \to {\Ch}_{\Q}^{\Out^{op}}$ for the left Kan extension along $q$. Concretely, $q_{!}$ quotients out by all inner conjugations. This process does not change the homology of complexes of the form $c(\Phi(X))$, since inner conjugations already act trivially on their homology.

Hence, the composite $T$ has the property that it turns rational geometric fixed point homotopy groups of an orthogonal spectrum into homology groups of the associated $\Out^{op}$-chain complex. In particular, it takes rational global equivalences to quasi-isomorphisms of $\Out^{op}$-chain complexes. In \cite{Wim17}, Wimmer shows that the induced functor on homotopy categories
\[ T:\textnormal{orthogonal spectra}[\Q-\textnormal{global equivalences}^{-1}]\to \mathcal{D}(\Q[{\Out}^{op}]-\textnormal{mod}) \]
is an exact equivalence of triangulated categories. We note that our results do not rely on this not yet published theorem, but only on the construction of $T$ and the properties above. Proposition \ref{prop:formal} also uses that $T$ preserves direct sums on the level of homotopy categories, which is a consequence of the fact that $\cC_*(-)$ takes wedges of cofibrant based spaces to direct sums of chain complexes, up to quasi-isomorphism.

Recall that by $\cC \widetilde{L}_n$ we denote the rational chains on the $\Out^{op}$-functor $G\mapsto \wL(G)_n/G$. Then the strongest version of the main result of this paper is the following:

\begin{Theorem} \label{theo:strong} There are quasi-isomorphisms of chain complexes of $\Q[\Out^{op}]$-modules
\[ T(Sp^n)\simeq \cC \widetilde{L}_n \]
for all $n\in \N$ that are compatible with the inclusions $Sp^n\to Sp^{n+1}$ and $\cC \widetilde{L}_n\to \cC \widetilde{L}_{n+1}$.
\end{Theorem}

Via taking homology and applying Proposition \ref{prop:equiv}, this implies Theorems \ref{thm:mackey} and \ref{thm:geom}.
\begin{Remark} \label{rem:gmorita} For fixed finite $G$, the category of rational $G$-spectra is Quillen-equivalent to the derived category of rational $G$-Mackey functors (which in turn is isomorphic to the product of the derived categories of rational $W_G(H)$-modules, see Remark \ref{rem:ggeom}). The category of rational $G$-Mackey functors is semisimple. Hence, a chain complex of such is determined up to quasi-isomorphism by its homology. Consequently, a rational $G$-spectrum is determined by its $G$-Mackey functor homotopy groups. It follows that over fixed $G$ the analogs of the statements of Theorems \ref{thm:geom} and \ref{theo:strong} are equivalent. However, since rational $\Out^{op}$-modules are not semisimple, these theorems are not equivalent globally. In Section \ref{sec:globalprop} we use Theorem \ref{theo:strong} to show that~$Sp^n_{\Q}$ is not a product of global Eilenberg-MacLane spectra unless $n$ is $1$ or~$\infty$.
\end{Remark}

\section{Proof of the equivalence} \label{sec:proof}
Now we come to the proof of Theorem \ref{theo:strong}. Let $|-|$ denote the geometric realization of a simplicial set. We show that there exists a transformation of spectrum-valued $\Epi^{op}$-functors \[ \widetilde{\alpha}:\Sigma^ {\infty} |\widetilde{L}(-)_n|\to \Phi (Sp^n) \]
which induces isomorphisms
\[ (H_*(\Sigma^{\infty}|\widetilde{L}(G)_n|,\Q))_G\cong H_*(\Phi^G (Sp^n),\Q) \]
for all finite groups $G$. Since rational homology is naturally isomorphic to rational stable homotopy, we see that this implies quasi-isomorphisms
\[ \cC\widetilde{L}_n=q_{!}(\cC_*(\widetilde{L}_n(-)))\simeq q_{!}(c(\Sigma^ {\infty}_+ |\widetilde{L}(-)_n|))\simeq q_{!}(c(\Phi (Sp^n)))=T(Sp^n) \]
of $\Out^{op}$-complexes and hence yields Theorem \ref{theo:strong}. Here, we made use of the properties of the functor $c$ that we described in the previous section.

More precisely, $\widetilde{\alpha}$ is a zig-zag, as we have to modify both $\wL(G)_n$ and $\Phi (Sp^n)$ to be able to construct an honest map. We now develop this zig-zag step by step, the whole construction is summarized in diagram \eqref{eq:zigzag} at the end of Section \ref{sec:naturality}. First we note that by adjunction, spectrum maps $\Sigma^{\infty} |\widetilde{L}(G)_n|\to \Phi(Sp^n)$ stand in bijection with maps of based spaces $|\wL(G)_n|\to (Sp^n(S^0))^G$. The target is just a discrete set of points, so there are no interesting maps on the point-set level. This can be resolved by stabilizing once, which we do via the following construction: The \emph{shift} $\shift X$ of an orthogonal spectrum $X$ is defined via $(\shift X)(V)=X(\R\oplus V)$, with structure maps the shifted ones of $X$. It allows a natural map $\lambda_X:S^1\wedge X \to \shift X$ given in level $V$ by the composite
\[ S^1\wedge X(V)\cong X(V)\wedge S^1\xr{\sigma_V^1} X(V\oplus \R)\xr{X(\tau_{\R,V})} X(\R\oplus V),\]
where $\tau_{\R,V}$ is the isometry $V\oplus \R\to \R\oplus V$ that swaps the two summands. The adjoint is a morphism $\widetilde{\lambda}_X:X\to \Omega \shift X$. Both $\lambda_X$ and $\widetilde{\lambda}_X$ induce isomorphisms on homotopy groups (cf. \cite[Proposition 3.1.25]{Sch18} for $G$ the trivial group).
Hence, instead of $\Phi (Sp^n)$ we can equivalently consider the $\Epi^{op}$-diagram $\Omega \shift \Phi (Sp^n)$. Morphisms $\Sigma^{\infty} |\wL(G)_n|\to \Omega \shift \Phi^G(Sp^n)$ now correspond to maps of based spaces $|\wL(G)_n|\to \Omega Sp^n(S^{\rho_G})^G$. This one copy of the regular representation turns out to be enough to define $\widetilde{\alpha}$ (though it only becomes a rational equivalence after further stabilization).

\subsection{Geometric idea} \label{sec:geomidea}

We start with the special case $n=|G|$ (and hence $L(G)_n=L(G)$) and first describe a map $\overline{\alpha}_1:|L(G)|\to \Omega (Sp^{|G|}(S^{\rho_G}))^G$ from the non-reduced subgroup lattice. The regular representation $\rho_G$ decomposes as $\R\oplus \overline{\rho}_G$, where $\overline{\rho}_G$ is the reduced regular representation of tuples that add up to $0$ and $\R$ denotes the trivial diagonal copy. To be explicit, we work with the splitting that sends an element $x=\sum x_g\cdot e_g$ to its trivial component $t(x)=(\sum x_g)\cdot (\frac{1}{|G|}\cdot \sum e_g)$ and its reduced component $r(x)=x-t(x)$. This decomposition also induces a map $\widetilde{\sigma}^G:(Sp^{|G|}(S^{\rrho_G}))^G\to \Omega (Sp^{|G|}(S^{\rho_G}))^G$, i.e., the induced map on fixed points of the adjoint structure map $\widetilde{\sigma}_{\rrho_G}^1$ of the orthogonal spectrum $Sp^{|G|}$. The map~$\overline{\alpha}_1$ that we construct is the composition of a map $\alpha_1:|L(G)|\to (Sp^{|G|}(S^{\rrho_G}))^G$ with $\widetilde{\sigma}^G$.

Given a subset $M\subseteq G$, we denote by $e_M\in \rho_G$ the element $\frac{1}{\sqrt{|M|}}\sum_{g\in M} e_g\in \rho_G$. Then $\alpha_1$ is defined by sending a chain of subgroups $H_0\leq \hdots \leq H_k$ of $G$ and $(t_0,\hdots,t_k)\in \Delta^k$ to the class
\[ [(r(\sum_{i=1}^k t_i\cdot e_{gH_i}))_{g\in G}]\in (Sp^{|G|}(S^{\rrho_G}))^G. \]

We shall explain this formula briefly: For each $g\in G$, the map
\[ (\alpha_1)_g:(H_0\leq \hdots \leq H_k;t_0,\hdots, t_k)\mapsto r(\sum_{i=1}^k t_i\cdot e_{gH_i})\]
defines an embedding of the subgroup lattice into $\rrho_G$. These different embeddings are permuted via the $G$-action, as an element $g'$ sends $(\alpha_1)_g$ to $(\alpha_1)_{g'g}$. So $\alpha_1$, the product of all $(\alpha_1)_g$, is $G$-fixed in $Sp^{|G|}(S^{\rrho_G})$.

\begin{Example} Elements in a symmetric product $Sp^n(X)$ of a space $X$ can be visualized as configurations in $X$ with labels in the natural numbers. The label on a point indicates how often it occurs in the tuple. We use this visualization to describe $\alpha_1$ in the case where $G$ is a cyclic group of order $3$ or $4$.

For $G=C_3$ the reduced regular representation is isomorphic to $\R^2$ with rotation by $120$ degrees. The image of the vertex $\{e\}$ in the subgroup lattice is the configuration of the three corners of an equilateral triangle with center $0$, each equipped with the label $1$. As one moves along the edge~$\{e\}\leq C_3$, these points move straight towards the center at the same speed. Finally, the vertex $C_3$ is mapped to the zero vector with label $3$. This map is described in Figure \ref{fig:c3}.

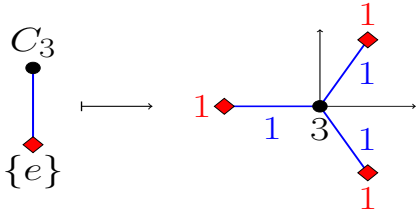
\begin{figure}
\huge
\centering
\resizebox{6cm}{3cm}{
\begin{tikzpicture}
\coordinate [label=above:$C_3$] (a) at (-6,1);
\coordinate [label=below:$\{e\}$] (b) at (-6,-1);
\draw[very thick,blue] (a) -- (b);
\node[draw,circle,inner sep=3pt,fill=black] at (a){};
\node[draw,diamond,inner sep=3pt,fill=red] at (b){};

\draw[|->,thick] (-5,0) -- (-3.5,0);

\coordinate (A) at (2,0);
\coordinate (B) at (0,2);
\coordinate [label={[label distance=.05cm,black]268:3}] (0) at (0,0);
\coordinate [label={[label distance=.1cm,red]above:1}] (X) at (60:2);
\coordinate [label={[red]left:1}] (Y) at (180:2);
\coordinate [label={[label distance=.1cm,red]below:1}] (Z) at (300:2);
\draw[->] (0)-- (A);
\draw[->] (0) -- (B);
%\draw (X)--(Y)--(Z)--cycle;
\draw[very thick,blue] (X) -- (0) node[midway, right] {1};
\draw[very thick,blue] (Y) -- (0) node[midway, below] {1};
\draw[very thick,blue] (Z) -- (0) node[midway, right] {1};
\node[draw,diamond,inner sep=3pt,fill=red] at (X){};
\node[draw,diamond,inner sep=3pt,fill=red] at (Y){};
\node[draw,diamond,inner sep=3pt,fill=red] at (Z){};
\node[draw,circle,inner sep=3pt,fill=black] at (0){};
\end{tikzpicture}
}
\normalsize
\caption{The map $\alpha_1$ for $G=C_3$}
\label{fig:c3}
\end{figure}

The reduced regular representation of $C_4$ is three-dimensional and permutes the corners of a regular tetrahedron. The image of the vertex $\{e\}$ under $\alpha_1$ is the configuration of these corners with label $1$, and the images of the other simplices are as depicted in Figure \ref{fig:c4}. Denoting a generator of $C_4$ by $t$, the rightmost corner corresponds to $r(e_1)$, the left one to $r(e_t)$, the upper one to $r(e_{t^2})$ and the lower one to~$r(e_{t^3})$.
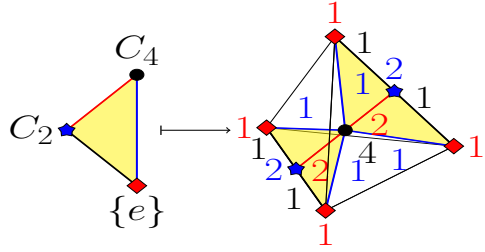
\begin{figure}
\huge
\centering 
\resizebox{7cm}{3.5cm}{
\begin{tikzpicture}
\node[coordinate,label=above:$C_4$] (a) at (-4,2.2,0) {};
\node[coordinate,label=below: $\{e\}$] (b) at (-4,-0.8,0) {};
\node[coordinate,label=left:$C_2$] (c) at (-5.5,0.7,0) {};
\draw[-, fill=yellow,, opacity=.5] (a)--(b)--(c)--cycle;
\draw[very thick,red] (a)--(c);
\draw[very thick,blue] (a)--(b);
\draw[very thick,black] (b)--(c);
\node[draw,circle,inner sep=3pt,fill=black] at (a){};
\node[draw,diamond,inner sep=3pt,fill=red] at (b){};
\node[draw,star,star points = 5,inner sep=2.5pt,fill=blue] at (c){};

\draw[|->,thick] (-3.5,0.7) -- (-2,0.7);

\coordinate [label={[red]right:1}] (A) at (2,-0.5,-2);
\coordinate [label={[red]left:1}] (B) at (-2,0,-2);
\coordinate [label={[label distance=0.1cm,red]above:1}] (C) at (1,4,2);
\coordinate [label={[label distance=0.1cm,red]below:1}] (D) at (0.8,-0.7,2);
\coordinate [label={[label distance=0.05cm,blue]above:2}] (E) at (1.5,1.75,0);
\coordinate [label={[label distance=0.05cm,black]275:4}] (F) at (0.45,0.7,0);
\coordinate [label={[blue]left:2}] (G) at (-0.6,-0.35,0);
\draw [thin] (A)--(D)--(B)--cycle;
\draw [thin] (B)--(D)--(C)--cycle;

\draw[-, fill=yellow, opacity=.5] (A)--(E)--(F)--cycle;
\draw[-, fill=yellow, opacity=.5] (C)--(E)--(F)--cycle;
\draw[-, fill=yellow, opacity=.5] (B)--(G)--(F)--cycle;
\draw[-, fill=yellow, opacity=.5] (D)--(G)--(F)--cycle;

\draw[very thick,black] (A) -- (E) node [midway,above] {1};
\draw[very thick,black] (C) -- (E) node [midway,above] {1};
\draw[very thick,black] (B) -- (G) node [midway,left] {1};
\draw[very thick,black] (D) -- (G) node [pos=0.3,left] {1};
\draw[very thick,blue] (A) -- (F) node [midway,below] {1};
\draw[very thick,blue] (B) -- (F) node [midway,above] {1};
\draw[very thick,blue] (C) -- (F) node [midway,right] {1};
\draw[very thick,blue] (D) -- (F) node [midway,right] {1};
\draw[very thick,red] (E) -- (F) node [pos=0.3,below] {2};
\draw[very thick,red] (G) -- (F) node [midway,below] {2};

\draw [thin] (A) --(D)--(C)--cycle;
\foreach \x in {A,B,C,D}{
\node[draw,diamond,inner sep=3pt,fill=red] at (\x){};
};
\node[draw,star,star points = 5,inner sep=2.5pt,fill=blue] at (E){};
\node[draw,star, star points = 5,,inner sep=2.5pt,fill=blue] at (G){};
\node[draw,circle,inner sep=3pt,fill=black] at (F){};
\end{tikzpicture}
}
\normalsize
\caption{The map $\alpha_1$ for $G=C_4$}
\label{fig:c4}
\end{figure}
\end{Example}

This is the basic geometric idea, but some adjustments are necessary in order to make it have all the properties and compatibilities that we need. One problem is that $\alpha_1$ does not yet factor through the reduced lattice $\wL(G)$, i.e., it does not send chains of subgroup inclusions which end in a proper subgroup of $G$ to the basepoint. This can be resolved as follows: Note that the full group $G$ is the only vertex that is sent to the $0$-vector. So if we choose a ball around $0$ of small enough radius and push everything that lies outside of it to $\infty$, the resulting map will send all proper subgroups $H$ of $G$ and the simplices connecting them to the basepoint.

To describe this in formulas, we let $p:S^{\rrho_G}\to S^{\rrho_G}$ be a map of the form $p(v)= \mu(|v|)\cdot v$, where $\mu$ is a fixed continuous self-map of $[0,\infty]$ that restricts to an orientation-preserving homeomorphism $[0,\frac{1}{\sqrt{2}}]\cong [0,\infty]$ and sends $[\frac{1}{\sqrt{2}},\infty]$ to $\infty$. In other words, $p$ collapses the hemisphere of vectors of length at least $\frac{1}{\sqrt{2}}$ to a point and identifies the resulting quotient with $S^{\rrho_G}$ again.
Furthermore, we let $q:(\rho_G-\{0\})\to S(\rho_G)$ denote the projection to the unit sphere and $\lr:(\rho_G-\{0\})\to \rrho_G$ the composite of $q$ and the retraction $r$ defined above. For every $g\in G$ we obtain a new map $(\alpha_2)_g:|L(G)|\to S^{\rrho_G}$ via the formula
\[ (H_0\leq \hdots \leq H_k;t_0,\hdots,t_k)\mapsto p(\lr(\sum_{i=0}^k t_i\cdot e_{gH_i})) \]
and again let
\[ \alpha_2:|L(G)|\to (Sp^{|G|}(S^{\rrho_G}))^G \]
be the tuple of all $(\alpha_2)_g$ for $g\in G$. In words, we have made two changes: We project each of the lattices inside $\rho_G$ to the unit sphere before passing to the reduced $\rrho_G$, and in the end we quotient out all vectors of length at least $\frac{1}{\sqrt{2}}$. This has the desired effect:
\begin{Lemma} \label{lem:reduced} The map $\alpha_2: |L(G)|\to (Sp^{|G|}(S^{\rrho_G}))^G$ factors through the reduced lattice~$|\wL(G)|$.
\end{Lemma}

\begin{proof}
It suffices to see that the square of the norm of 
\[ q(\sum_{i=0}^k (t_i\cdot e_{gH_i}))+ t\cdot e_G \]
is at least $1/2$ for any chain $H_0\leq \hdots \leq H_k$, $(t_0,\hdots,t_k)\in \Delta^k$ and $t\in \R$, provided that~$H_k$ is a proper subgroup of $G$.
Dividing $\rho_G$ into the span of the basis elements of the form $e_{gh_k}$ with $h_k\in H_k$ and the span of the other basis elements, we see that this square is given by
\[ |q(\sum_{i=0}^k (t_i\cdot e_{gH_i}))+t\cdot \sqrt{\frac{|H_k|}{|G|}} \cdot e_{gH_k}|^2 + |t\cdot \sqrt{\frac{|G|-|H_k|}{|G|}}\cdot e_{(G-gH_k)}|^2. \]
Using that $|e_{gH_k}|=1=|e_{(G-gH_k)}|$ and applying the triangle inequality yields that this square is at least as large as
\[ \left(|q(\sum_{i=0}^k (t_i\cdot e_{gH_i}))|-|t|\cdot \sqrt{\frac{|H_k|}{|G|}}\right)^2 + \left(|t|\cdot \sqrt{\frac{|G|-|H_k|}{|G|}}\right)^2. \]
Since $q(-)$ by definition always has norm $1$ and $|G|-|H_k|$ is at least $|H_k|$, we obtain the lower bound
\[ \left(1-|t|\cdot \sqrt{\frac{|H_k|}{|G|}}\right)^2 + \left(|t|\cdot \sqrt{\frac{|H_k|}{|G|}}\right)^2.\]
The minimum of this quadratic function equals $\frac{1}{2}$, which proves the claim.
\end{proof}

\begin{Example} \label{exa:alpha2} The effect of $\alpha_2$ is depicted in Figure \ref{fig:alpha2} for $G=C_3$. The first image illustrates the area of $\rrho_{C_3}$ that is quotiented out and the second the resulting map to $Sp^3(S^{\rrho_{C_3}})$.
\end{Example}
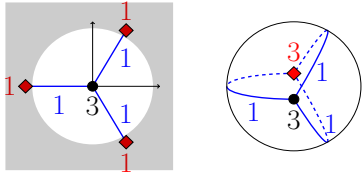
\begin{figure}
\huge
\centering
\resizebox{5cm}{2.5cm}{\begin{tikzpicture}

\coordinate (A) at (2,0);
\coordinate (B) at (0,2);
\coordinate [label={[label distance=0.05cm,black]below:3}] (0) at (0,0);
\coordinate [label={[label distance=0.05cm,red]above:1}] (X) at (60:2);
\coordinate [label={[red]left:1}] (Y) at (180:2);
\coordinate [label={[label distance=0.1cm,red]below:1}] (Z) at (300:2);
\draw[->] (0)-- (A);
\draw[->] (0) -- (B);
\draw[very thick,blue] (X) -- (0) node[midway, right] {1};
\draw[very thick,blue] (Y) -- (0) node[midway, below] {1};
\draw[very thick,blue] (Z) -- (0) node[midway, right] {1};
\node[draw,diamond,inner sep=3pt,fill=red] at (X){};
\node[draw,diamond,inner sep=3pt,fill=red] at (Y){};
\node[draw,diamond,inner sep=3pt,fill=red] at (Z){};
\node[draw,circle,inner sep=3pt,fill=black] at (0){};

\fill[black,opacity=.2,even odd rule]
    (-2.6,2.6) rectangle (2.4,-2.6)
    (0,0) circle[radius=1.8];

\draw[very thick,blue,rotate around={60:(6,0)}] (6,0) [partial ellipse=258:360:2cm and 0.2cm];
\draw[very thick,blue,rotate around={60:(6,0)},dashed] (6,0) [partial ellipse=360:440:2cm and 0.2cm];
\draw[very thick,blue,rotate around={120:(6,0)}] (6,0) [partial ellipse=100:180:2cm and 0.2cm];
\draw[very thick,blue,rotate around={120:(6,0)},dashed] (6,0) [partial ellipse=180:280:2cm and 0.2cm];
\draw[very thick,blue,dashed] (6,0) [partial ellipse=90:180:2cm and 0.4cm];
\draw[very thick,blue] (6,0) [partial ellipse=180:270:2cm and 0.4cm];

\draw (6,0) circle (2cm);
\coordinate [label={[label distance = 0.1cm,red]above:3}] (02) at (6,0.4);
\coordinate [label={[label distance = 0.1cm,black]below:3}] (x) at (6,-0.4);
\node[draw,diamond,inner sep=3pt,fill=red] at (6,0.4){};
\node[draw,circle,inner sep=3pt,fill=black] at (6,-0.4){};

\coordinate [label={[blue]below:1}] (x1) at (4.8,-0.25);
\coordinate [label={[blue]below:1}] (x1) at (6.85,1);
\coordinate [label={[blue]below:1}] (x1) at (7.1,-0.6);

\end{tikzpicture} }
\normalsize
\caption{The map $\alpha_2$ for $G=C_3$}
\label{fig:alpha2} 
\end{figure}
However, there is a problem that is more complicated to resolve: We need the restriction to the subcomplex $\widetilde{L}(G)_n\subseteq \widetilde{L}(G)$ to take image in $(Sp^n(S^{\rrho_G}))^G$.
This is simply not the case for $\alpha_2$, as one already sees in Figure \ref{fig:alpha2}: The image of the $0$-chain with value $C_3$ (which lies in $\wL(G)_1$) is the tuple $[(0,0,0)]$. It has three non-basepoint components and hence does not lie in any smaller symmetric product. The idea to rectify this is to use that $[(0,0,0)]$ is stably the same as `three times the element $[(0,*,*)]$', which does lie in the image of $Sp^1$.

We now make this precise and more generally let $H_0\leq \hdots \leq H_k$ be a chain of total index $n$ and ending in $H_k=G$. We write
\[ (\alpha_2)_g(\{H_i\},-):\Delta^k\to S^{\rrho_G} \]
for the restriction of $(\alpha_2)_g$ to the $k$-simplex corresponding to this chain of subgroup inclusions, and
\[\alpha_2(\{H_i\},-):\Delta^k\to (Sp^{|G|}(S^{\rrho_G}))^G\]
for the analogous restriction of $\alpha_2$. These have the following properties:
\begin{itemize} \label{items:notation}
\item  Each $(\alpha_2)_g(\{H_i\},-)$ only depends on the coset $gH_0$, since multiplication with $h_0\in H_0$ leaves all $e_{H_i}$ fixed. Hence there are only $n$ different components in $\alpha_2(\{H_i\},-)$, each repeated $|H_0|$ times.
\item Let $g_1,\hdots,g_n$ be a system of coset representatives of $G/H_0$. Then the tuple 
\[ (\alpha_2)_{G/H_0}(\{H_i\},-)\defeq [((\alpha_2)_{g_j}(\{H_i\},-))_{j=1,\hdots,n}] \]
defines a  map $\Delta^k\to (Sp^{n}(S^{\rrho_G}))^G$. This map does not depend on the choice of the $g_i$.
\end{itemize}
In other words, $\alpha_2(\{H_i\},-)$ factors through the diagonal 
\[ (Sp^{n}(S^{\rrho_G}))^G\xr{\Delta} (Sp^{|G|}(S^{\rrho_G}))^G \]
that repeats each entry $|H_0|$ times, while we want it to factor through the standard inclusion $i_{n}^{|G|}$. In $(Sp^{|G|}(S^{\rrho_G}))^G$ there is no direct way to pass between $\Delta$ and $i_{n}^{|G|}$, but there is after stabilizing once, i.e., after postcomposing with 
\[ \widetilde{\sigma}^G:(Sp^{|G|}(S^{\rrho_G}))^G\to \Omega (Sp^{|G|}(S^{\rho_G}))^G. \]
To see this, we consider the following modified construction of the diagonal: We assume given $|H_0|$ closed subintervals $[a_i,b_i]$ of $[-\infty,\infty]$ and for each of these let $c(a_i,b_i)$ denote the self-map of $S^1$ which collapses everything outside $(a_i,b_i)$ to the basepoint and identifies $[a_i,b_i]$ with $[-\infty,\infty]$ in some fixed orientation-preserving way. To each such data one can associate a map 
\[ \Delta_{\{[a_i,b_i]\}}:\Omega (Sp^{n}(S^{\rho_G}))^G\to \Omega (Sp^{|G|}(S^{\rho_G}))^G \]
by sending $\varphi\in \Omega (Sp^{n}(S^{\rho_G}))^G$ to
\[ [(\varphi\circ c(a_i,b_i))_{1\leq i \leq |H_0|}]. \]
So, instead of repeating it $|H_0|$ times, $\varphi$ is precomposed with every $c(a_i,b_i)$. If all of the $[a_i,b_i]$ are equal to $[-\infty,\infty]$ and one takes the identity of $[-\infty,\infty]$ for the identification, this construction gives back the usual diagonal~$\Omega(\Delta)$. If on the other hand the interiors of the $[a_i,b_i]$ are pairwise disjoint, the map $\Delta_{\{[a_i,b_i]\}}$ factors through the standard inclusion $\Omega(i_{|G/H|}^{|G|})$. Indeed, in this case there is at most one $i$ such that $c(a_i,b_i)(t)$ is not the basepoint, for any fixed $t\in \R$. So each loop $\Delta_{\{[a_i,b_i]\}}(\varphi)$ has at most $n$ non-trivial components at every $t$. This means that it has at most $n$ non-trivial components globally and hence lies in the image of $\Omega(i_{n}^{|G|})$, since we can always move the non-trivial components to the first $n$ entries. The self-map of $\Omega Sp^{n}(S^{\rho_G})$ obtained this way can also be described differently:
It is given by precomposition with the self-map of $S^1$ that collapses everything outside the open intervals $(a_i,b_i)$ to the basepoint and identifies each $[a_i,b_i]/(a_i\sim b_i)$ with $S^1$. In particular, the homotopy class of $\Delta_{\{[a_i,b_i]\}}(\varphi)$ is the $|H_0|$-fold sum of $\varphi$ with itself.

Any choice of homotopies from the $c(a_i,b_i)$ to the identity of $S^1$ induces a homotopy between $\Omega (\Delta)$ and $\Delta_{\{[a_i,b_i]\}}$. So we see that, up to reparametrization of loops and in particular up to homotopy, $\widetilde{\sigma}^G\circ \alpha_2$ does map the simplex associated to $H_0\leq \hdots \leq H_k$ to the image of $\Omega Sp^{n}(S^{\rho_G})$ under $\Omega(i_{n}^{|G|})$. To turn this into honest maps from $\widetilde{L}(G)_n$ to $\Omega(Sp^n(S^{\rho_G}))^G$, we need to make choices of reparametrizations that are coherent for all chains $H_0\leq \hdots \leq H_k$, all $n\in \N$ and all finite groups $G$.  We deal with this by defining a modification of the subgroup lattice that contains a contractible choice of intervals as part of the data.

\subsection{Fattening of the lattice}
For this it turns out to be more convenient to work with subintervals of $[0,1]$ instead of $[-\infty,\infty]$. Whenever we need to switch between the two, we use the homeomorphism that maps $t\in [0,1]$ to $\frac{2t-1}{t(1-t)}\in [-\infty,\infty]$.

Let $J$ denote the space of closed subintervals $[a,b]$ of $[0,1]$ (with $a<b$), topologized as a subspace of $[0,1]\times [0,1]$. We let $L_f(G)$ denote the following topological category: The object space is given by 
\[ \bigsqcup_{H\leq G} Sp^{|H|}(J),\]
i.e., subgroups $H$ of $G$ together with an unordered $|H|$-tuple of subintervals of $[0,1]$. The morphism space from a component $Sp^{|H|}(J)$ to $Sp^{|K|}(J)$ is empty if $H$ is not contained in $K$ and is otherwise given by $Sp^{|H|}(J)$ again.
In this case the target map 
\[ Sp^{|H|}(J)\to Sp^{|K|}(J) \]
is the diagonal which repeats each subinterval $[K:H]$ times, and the source map is the identity. Hence, a $k$-simplex in the topological nerve of $L_f(G)$ - which we also denote by $L_f(G)_n$ - is given by a chain of subgroups $H_0\leq \hdots \leq H_k$ together with $|H_0|$ many subintervals $[a_i,b_i]$ of $[0,1]$. We filter this nerve by saying that such a $k$-simplex lies in $L_f(G)_n$ if
\begin{enumerate}[(i)]
	\item the total index $[H_k:H_0]$ is at most $n$, and
	\item the intervals $([a_i,b_i])$ have at most $\frac{n}{[H_k:H_0]}$-fold intersections, i.e., every $t\in (0,1)$ lies in the interior of at most $\frac{n}{[H_k:H_0]}$-many $[a_i,b_i]$.
\end{enumerate}

There is an obvious forgetful functor $\mu:L_f(G)\to L(G)$ to the usual subgroup lattice of $G$, whose nerve maps $L_f(G)_n$ into $L(G)_n$.
\begin{Lemma} The maps $\mu:L_f(G)_n\to L(G)_n$ induce homotopy equivalences on geometric realizations.
\end{Lemma}
\begin{proof} Given $l,m\in \N$, let $J^{l}_{m}$ denote the subspace of $Sp^l(J)$ of tuples of intervals with at most $m$-fold intersections. Then the space of $k$-simplices of $L_f(G)_n$ is given by the disjoint union
\begin{equation} \label{eq:J} \bigsqcup_{H_0\leq \hdots \leq H_k; [H_k:H_0]\leq n} J^{|H_0|}_{\lfloor \frac{n}{[H_k:H_0]} \rfloor},
\end{equation}
where $\lfloor \frac{n}{[H_k:H_0]} \rfloor$ denotes the largest integer smaller than or equal to $\frac{n}{[H_k:H_0]}$. We first claim that each $J^l_{m}$ is contractible and hence $L_f(G)_n\to L(G)_n$ forms a degreewise homotopy equivalence. Since we have quotiented out by the symmetric group action, every element of $J^l_{m}$ has a unique representative $([a_1,b_1],\hdots,[a_l,b_l])$ for which the $(a_i,b_i)$ are lexicographically ordered, i.e., $a_i\leq a_{i+1}$ and if $a_i=a_{i+1}$ then $b_i\leq b_{i+1}$. In fact, $J^l_{m}$ is homeomorphic to the space of such ordered tuples with at most $m$-fold intersections. But this space is star-shaped, it can be linearly contracted onto the tuple $([0,\frac{1}{l}],[\frac{1}{l},\frac{2}{l}],\hdots,[\frac{l-1}{l},1])$. This proves the claim.

Hence it suffices to see that both $L_f(G)_n$ and $L(G)_n$ are Reedy cofibrant simplicial spaces (cf. Hischhorn \cite[Chapter 15]{Hir03}) with respect to the Str\o m model structure on topological spaces (\cite{Str72}). For $L(G)_n$ this is clear, since it is a discrete simplicial space. For $L_f(G)_n$, the $k$-th latching map is the inclusion of those components in the disjoint union \eqref{eq:J} above which are associated to chains $H_0\leq \hdots \leq H_k$ for which at least one containment is not proper. Every topological space is Str\o m cofibrant, so this inclusion is a Str\o m cofibration, which finishes the proof.
\end{proof}

Hence, the $L_f(G)_n$ indeed form a fattening of the $L(G)_n$, but we still need to explain their functoriality in surjective group homomorphisms. For this we let $\psi:G\twoheadrightarrow K$ be a surjection and denote by $k$ the order of the kernel. Then we define 
\[ \psi^*:L_f(K)\to L_f(G) \]
to send a subgroup $L$ of $K$ to $\psi^{-1}(L)$, and the associated collection of intervals $([a_i,b_i])$ to 
\[ (r^k_j([a_i,b_i]))_{i=1,\hdots,l;j=0,\hdots,k-1}, \]
where $r^k_j:[0,1]\hookrightarrow [0,1]$ is the unique oriented affine embedding with image $[\frac{j}{k},\frac{j+1}{k}]$. In other words, $\psi^*$ splits $[0,1]$ into $k$ parts of equal size and copies each $[a_i,b_i]$ into every one of them, yielding $|L|\cdot k=|\varphi^{-1}(L)|$ subintervals, as needed.
This definition turns $L_f(-)$ into a functor from $\Epi^{op}$ to topological categories. After applying the nerve, it restricts to functors $L_f(-)_n$ from $\Epi^{op}$ to simplicial spaces and hence after geometric realization to functors $|L_f(-)_n|:\Epi^{op}\to T$. 
The forgetful functor $|L_f(-)_n|\to |L(-)_n|$ is natural for this $\Epi^{op}$-functoriality and a homotopy equivalence for all finite groups~$G$.

Finally, we again define a reduced version $|\widetilde{L}_f(G)_n|$ by quotienting out all simplices associated to chains that do not end in the full group $G$. The $|\widetilde{L}_f(G)_n|$ again assemble to a functor $\Epi^{op}\to T$ and the forgetful map $|\widetilde{L}_f(-)_n|\to |\widetilde{L}(-)_n|$ defines a natural levelwise based homotopy equivalence.

\subsection{Definition} \label{sec:construction}
Given an interval $[a,b]$ of $[0,1]$, we from now on let $c(a,b)$ denote the self-map of $[0,1]/\{0,1\}$ (or $S^1=\R\cup \{\infty \}$, using the fixed homeomorphism between the two) obtained by collapsing everything outside $(a,b)$ to a point and using the identification $[a,b]\cong [0,1]$ that sends $x$ to~$\frac{x-a}{b-a}$.

We are now ready to define the map 
\[ \alpha:|L_f(G)| \to \Omega (Sp^{|G|}(S^{\rho_G}))^G \]
by sending a simplex associated to a chain $H_0\leq \hdots \leq H_k$ together with intervals $[a_1,b_1],\hdots, [a_{|H_0|},b_{|H_0|}]$ to the composite
\begin{equation} \label{eq:alphadelta} \Delta^k \xr{(\alpha_2)_{G/H_0}(\{H_i\},-)} (Sp^{|G/H_0|}(S^{\rrho_G}))^G\xr{\widetilde{\sigma}^G} \Omega (Sp^{|G/H_0|}(S^{\rho_G}))^G\xr{\Delta_{\{[a_i,b_i]\}}} \Omega (Sp^{|G|}(S^{\rho_G}))^G. \end{equation}
The maps $(\alpha_2)_{G/H_0}(\{H_i\},-)$ and $\Delta_{\{[a_i,b_i]\}}$ are explained in the last part of Section \ref{sec:geomidea} (for which we use the specific $c(a_i,b_i)$ defined above to construct the diagonal). If the $[a_j,b_j]$ are all equal to $[0,1]$, we get back $\widetilde{\sigma}^G\circ \alpha_2(\{H_i\},-)$.

In order for the maps \eqref{eq:alphadelta} to glue to a map from the geometric realization, we need to check that they are still compatible with the simplicial structure maps. This is a consequence of the fact that the $(\alpha_2)_g(\{H_i\},-)$ have this compatibility, except for the boundary $d_0$, since it changes the smallest subgroup $H_0$. We recall that $d_0^*$ of a tuple $(H_0\leq \hdots \leq H_k,\{[a_i,b_i]\})$ is the chain $H_1\leq \hdots \leq H_k$ together with the intervals $\Delta(\{[a_i,b_i]\})$, i.e., each interval repeated $|H_1/H_0|$ times. So the compatibility for $d_0$ follows from the commutativity of the diagram:
\[ \xymatrix{ \Delta^k \ar[r]^-{(\alpha_2)_{G/H_0}} & (Sp^{|G/H_0|}(S^{\rrho_G}))^G\ar[r]^{\widetilde{\sigma}^G} & \Omega (Sp^{|G/H_0|}(S^{\rho_G}))^G\ar[rr]^{\Delta_{\{[a_i,b_i]\}}} && \Omega (Sp^{|G|}(S^{\rho_G}))^G \\
\Delta^{k-1} \ar[u]_{d_0}\ar[r]^-{(\alpha_2)_{G/H_1}} & (Sp^{|G/H_1|}(S^{\rrho_G}))^G\ar[r]^{\widetilde{\sigma}^G}  \ar[u]^{\Delta} & \Omega (Sp^{|G/H_1|}(S^{\rho_G}))^G\ar[rr]^{\Delta_{\Delta(\{[a_i,b_i]\})}} \ar[u]^{\Omega (\Delta)} && \Omega (Sp^{|G|}(S^{\rho_G}))^G \ar[u]^{=}}
\]

Since $(\alpha_2)_{G/H_0}(\{H_i\},-)$ sends all chains that do not end in the full group $G$ to the basepoint (Lemma \ref{lem:reduced}), it follows that $\alpha$ again factors through the reduced fat lattice, yielding a map $\walph:|\widetilde{L}_f(G)|\to \Omega (Sp^{|G|}(S^{\rho_G}))^G$. In addition, we now have:

\begin{Prop} The restriction of $\walph$ to $|\widetilde{L}_f(G)_n|$ factors through 
\[ \Omega(i_n^{|G|}):\Omega (Sp^n(S^{\rho_G}))^G\hookrightarrow \Omega (Sp^{|G|}(S^{\rho_G}))^G. \]
\end{Prop}

\begin{proof}
Let $H_0\leq \dots \leq H_k$ be a chain with $H_k=G$ and $[G:H_0]\leq n$, together with $|H_0|$ intervals $[a_i,b_i]$ with at most $\frac{n}{[G:H_0]}$-fold intersections. Then at any point $t\in (0,1)$, at most $\frac{n}{[G:H_0]}$-many values $c(a_i,b_i)(t)$ are not equal to the basepoint. Each one of them appears exactly $[G:H_0]$ times in the definition of $\Delta_{\{[a_i,b_i]\}}$, so it follows that the diagonal 
\[ \Omega (Sp^{|G/H_0|}(S^{\rho_G}))^G\xr{\Delta_{\{[a_i,b_i]\}}} \Omega (Sp^{|G|}(S^{\rho_G}))^G \]
factors through $\Omega (i_n^{|G|})$, which proves the claim.
\end{proof}

\subsection{Naturality} \label{sec:naturality}
By adjunction, we obtain maps 
$ \walph:\Sigma^{\infty} |\widetilde{L}_f(G)_n|\to \Omega \shift \Phi^G(Sp^n)$,
compatible with the respective inclusions from $n$ to $n+1$. We now check their naturality with respect to surjective group homomorphisms $\psi:G\twoheadrightarrow K$.

Let $M$ be a subset of $K$. Then the linear isometry $i_{\psi}:\psi^*(\rho_K)\hookrightarrow \rho_G$ (defined in Section \ref{sec:globalchain} to describe the $\Epi^{op}$-functoriality of geometric fixed points) sends the element $e_M$ to $e_{\psi^{-1}(M)}$. This implies that for every chain of subgroups $H_0\leq \hdots \leq H_k$ of $K$, the composite 
\[ \Delta^k \xr{(\alpha_2)_{K/H_0}(\{H_i\},-)} (Sp^{|K/H_0|}(S^{\rrho_K}))^K \xr{(i_{\psi})_*} (Sp^{|G/\psi^{-1}(H_0)|}(S^{\rrho_G}))^G \]
equals the map $(\alpha_2)_{G/\psi^{-1}(H_0)}(\{\psi^{-1}(H_i)\},-)$. To compute the effect of $(\Omega \shift \Phi^{\psi})\circ \walph$ on the tuple $(H_0\leq \hdots\leq H_k,\{[a_i,b_i]\})$, we then have to postcompose with the diagonal~$\Delta_{\{[a_i,b_i]\}}$. On the other hand, in order to compute the effect of $\walph$ on the tuple $\psi^*(H_0\leq \hdots\leq H_k,\{[a_i,b_i]\})$ we have to postcompose with the diagonal $\Delta_{\{[r_k^l(a_i),r_k^l(b_i)]\}}$, where $k$ denotes the order of the kernel of $\psi$. The diagonal $\Delta_{\{[r_k^l(a_i),r_k^l(b_i)]\}}$ can be written as the composite
\[ \Delta_{\{[\frac{l-1}{k},\frac{l}{k}]\}}\circ \Delta_{\{[a_i,b_i]\}}. \]
A priori, this composite takes image in $\Omega (Sp^{|G|}(S^{\rho_G}))^G$, but since the intervals $(\frac{l-1}{k},\frac{l}{k})$ are pairwise disjoint, the diagonal $\Delta_{\{[\frac{l-1}{k},\frac{l}{k}]\}}$ factors as
\[ \xymatrix{\Omega (Sp^{|K|}(S^{\rho_G}))^G \ar[rr]^{\Delta_{\{[\frac{l-1}{k},\frac{l}{k}]\}}} \ar[d]_{l_k}&& \Omega (Sp^{|G|}(S^{\rho_G}))^G \\
\Omega (Sp^{|K|}(S^{\rho_G}))^G \ar[rru]_{\Omega (i_{|K|}^{|G|})}}
\]
where $l_k$ is the selfmap of $[0,1]/\{0,1\}$ which takes $x$ to $kx-\lfloor kx \rfloor$.

So, in summary, $\walph$ is natural for a different $\Epi^{op}$-functoriality on $\Omega \shift \Phi(Sp^n)$ which sends a surjection $\psi:G\twoheadrightarrow K$ to
\[ \Omega \shift \Phi^K(Sp^n)\xr{\Omega \shift^{\psi}} \Omega \shift \Phi^G(Sp^n) \xr{l_k^*} \Omega \shift \Phi^G(Sp^n), \]
where $l_k^*$ is the self-map of $\Omega \shift \Phi^G(Sp^n)$ which precomposes each loop with $l_k$. This difference in the functoriality can be corrected: The self-maps 
\[ l_{|G|}^*:\Omega \shift \Phi^G(Sp^n)\to \Omega \shift \Phi^G(Sp^n) \] assemble to a transformation from the usual $\Epi^{op}$-functor $\Omega \shift \Phi$ to this twisted one. Since $l_{|G|}$ induces multiplication by $|G|$ on homotopy, the transformation is a rational equivalence.
\begin{Remark} This `defect' of $\walph$ can be explained: Recall from Remark \ref{rem:fraction} that on $\pi_0$ we should be sending a vertex $H$ of the subgroup lattice to the element $\frac{1}{[G:H]} \cdot [G/H]$ in the rationalized Burnside ring. But this is impossible, since $\alpha$ is geometrically defined to land in the not yet rationalized spectrum $\Omega \shift \Phi(Sp^n)$. Instead it sends $H$ to $|H|\cdot [G/H]$, which needs to be corrected by dividing by $|G|$ afterwards.
\end{Remark}
So we finally obtain a zig-zag of natural transformations of $\Epi^{op}$-functors
\begin{equation} \label{eq:zigzag} \Sigma^{\infty} |\wL(-)_n| \xleftarrow{\simeq} \Sigma^{\infty} |\widetilde{L}_f(-)_n|\xr{\walph} (\Omega \shift \Phi(Sp^n))^{\text{twisted}} \xleftarrow{\simeq _{\Q}} \Omega \shift \Phi(Sp^n) \xleftarrow{\simeq} \Phi(Sp^n). \end{equation}
In order to prove Theorem \ref{theo:strong} it now remains to show that for all $n$ and $G$ the map $\walph$ induces an isomorphism  \[ (H_*(|\widetilde{L}_f(G)_n|,\Q))_G\xr{\cong} H_*(\Omega \shift \Phi^G(Sp^n),\Q). \]
Via an induction on $n$ and the five-lemma, this in turn can be reduced to showing that $\walph$ induces isomorphisms
\begin{equation} \label{eq:quotients} (H_*(|\widetilde{L}_f(G)_n/\widetilde{L}_f(G)_{n-1}|,\Q))_G\xr{\cong} H_*(\Omega \shift \Phi^G(Sp^n/Sp^{n-1}),\Q), \end{equation}
and this is what we will do.

\subsection{Rational splitting on subquotients} \label{sec:split}
Our first aim is to show:
\begin{Prop} \label{prop:splitinj} The map (\ref{eq:quotients}) above is split injective for all $n\in \N$ and finite groups~$G$.
\end{Prop}
We produce this splitting geometrically. Let $[(x_1,\hdots,x_n)]\in Sp^n(S^{k\cdot \rho_G})$ be a $G$-fixed point. Then the subset $\{x_1,\hdots,x_n\}\subseteq S^{k\cdot \rho_G}$ is closed under the $G$-action. Let $C^n$ be the subspectrum of $\Phi^G(Sp^n)$ consisting of those $G$-fixed points for which this $G$-set is not transitive or contains less than $n$ elements. In particular, $C^n$ contains $\Phi(Sp^{n-1})$ and so we can consider the composite
\[ \Sigma^{\infty}|\wL_f(G)_n/\wL_f(G)_{n-1}| \xr{\walph} \Omega \shift \Phi^G(Sp^n/Sp^{n-1})\to \Omega \shift (\Phi^G(Sp^n)/C^n). \]
Our aim is to show that this composite induces an isomorphism \[ (H_*(|\widetilde{L}_f(G)_n/\widetilde{L}_f(G)_{n-1}|,\Q))_G\xr{\cong} H_*(\Omega \shift (\Phi^G(Sp^n)/C^n),\Q), \]
which proves Proposition \ref{prop:splitinj}.

Every non-basepoint element of $(\Phi^G(Sp^n)/C^n)_k$ is determined by any of its components $x_i\in S^{k\cdot \rho_G}$, and the isotropy of such a point is necessarily an index $n$ subgroup of~$G$. Given an index $n$ subgroup $H$, we let $S(k,H)$ denote the space 
\[ (S^{k\cdot \rho_G})^H/\left(\colim_{H\lneq K\leq G} (S^{k\cdot \rho_G})^K\right), \]
i.e., the $H$-fixed points of the $G$-space $S^{k\cdot \rho_G}$ modulo all fixed points of larger subgroups. The spaces $S(k,H)$ assemble to an orthogonal spectrum $S(H)$, with structure map sending $(x\wedge t)\in (S^{k\cdot \rho_G})^H\wedge S^1$ to the class of the element $(x+t\cdot e_G) \in S^{(k+1)\cdot \rho_G}$. There are morphisms
\[ S(H)\to \Phi^G(Sp^n)/C^n \]
sending $x$ to $[(x,g_1\cdot x,\hdots,g_{n-1}\cdot x)_{[g_i]\in G/H}]$. As we just argued, every element $x=[(x_1,\hdots,x_n)]\in S^{k\cdot \rho_G}$ lies in the image of one of these, by choosing a component $x_i$. Moreover, the choice of a different component amounts to multiplying with some element $g\in G$, since the $G$-set $\{x_1,\hdots,x_n\}$ 
is assumed to be transitive. So we find that there is an 
isomorphism of orthogonal spectra
\begin{equation} \label{eq:split1} \Phi^G(Sp^n)/C^n\cong \left(\bigvee_{H\leq G,[G:H]=n} S(H)\right)/G,
\end{equation}
where the modded out $G$-action is given by translation, it sends $x\in (S^{k\cdot \rho_G})^H$ to $g\cdot x\in (S^{k\cdot \rho_G})^{gHg^{-1}}$.

\begin{Remark} This translation action should not be confused with the conjugation action on $\Phi^G(Sp^n)/C^n$ that comes from the $\Epi^{op}$-functoriality of geometric fixed points (as described in Section \ref{sec:globalchain}). The conjugation action sends an element $x\in (S^{k\cdot \rho_G})^H$ to $g\cdot x\cdot g^{-1}\in (S^{k\cdot \rho_G})^{gHg^{-1}}$, using that $\rho_G$ is both a left and a right module over~$G$. Together the conjugation action and the translation action assemble to an action of the semi-direct product $G\ltimes G$ on $\bigvee S(H)$. Under the isomorphism \eqref{eq:split1} above, the conjugation action on $\Phi^G(Sp^n)/C^n$ is the induced one on translation orbits. To make clear which action we are talking about, we write $-/^tG$ for the quotient by the translation action and $-/^cG$ for the one by the conjugation action.
\end{Remark}
We again use the decomposition $\rho_G\cong \R\oplus \overline{\rho}_G$ to rewrite each $(S^{k\cdot \rho_G})^H$ as $S^k\wedge (S^{k\cdot \overline{\rho}_G})^H$. This induces a decomposition $S(k,H)\cong S^k\wedge \overline{S}(k,H)$ with $\overline{S}(k,H)$ defined similarly to $S(k,H)$, replacing each $\rho_G$ by~$\overline{\rho}_G$. Through this identification, the structure map of $S(H)$ becomes the smash product of the associativity isomorphism $S^k\wedge S^1\cong S^{k+1}$ and the closed inclusion $\overline{S}(k,H)\hookrightarrow \overline{S}(k+1,H)$. Hence, the inclusions $\overline{S}(k,H)\hookrightarrow \overline{S}(\infty,H)$ assemble to a morphism of spectra
\begin{equation} \label{eq:sh} S(H)\to \Sigma^{\infty}\overline{S}(\infty,H). \end{equation}
By cofinality, this morphism induces an isomorphism on homotopy groups and is hence a stable equivalence.

We now turn to $\wL_f(G)_n/\wL_f(G)_{n-1}$, or rather its non-fat version $\wL(G)_n/\wL(G)_{n-1}$. Non-basepoint simplices in this quotient are given by chains of subgroup inclusions which go from an index $n$ subgroup $H$ to $G$. So we find that there is an isomorphism
\begin{equation} \label{eq:split2} \wL(G)_n/\wL(G)_{n-1}\cong \bigvee_{\substack{H\leq G \\ [G:H]=n}}\left(L(G)^{[H,G]}/\{\text{non-max. chains}\}\right), \end{equation}
where we write $L(G)^{[H,G]}$ for the poset of subgroups of $G$ which contain $H$ and $\{\text{non-max. chains}\}$ denotes the subcomplex of chains of subgroups that do not start in $H$ or do not end in $G$. The homotopy equivalence $L_f(G)_n^{[H,G]}\xr{\simeq} L(G)^{[H,G]}$ is split by the functor which sends a subgroup $K$ between $H$ and $G$ to itself together with the intervals $[0,\frac{1}{|H|}],[\frac{1}{|H|},\frac{2}{|H|}],\hdots,[\frac{|H|-1}{|H|},1]$, each repeated $[K:H]$ times.
Taking the wedge over these, we obtain a homotopy equivalence \[ |\widetilde{L}(G)_n/\widetilde{L}(G)_{n-1}|\cong \bigvee_{\substack{H\leq G \\ [G:H]=n}}\left(|L(G)^{[H,G]}|/\{\text{non-max. chains}\}\right)\xr{\simeq} |\wL_f(G)_n/\wL_f(G)_{n-1}|.\]
Composing this equivalence with the morphism 
\[ \Sigma^{\infty}|\wL_f(G)_n/\wL_f(G)_{n-1}|\to \Omega \shift (\Phi(Sp^n)/C^n)\]
and using \eqref{eq:split1}, \eqref{eq:sh} and \eqref{eq:split2}, we obtain a morphism
\begin{equation} \label{eq:comp} \bigvee_{\substack{H\leq G \\ [G:H]=n}}\Sigma^{\infty} \left(L(G)^{[H,G]}/\{\text{non-max. chains}\}\right)\to \Omega \shift ((\bigvee_{\substack{H\leq G \\ [G:H]=n}}\Sigma^{\infty} \overline{S}(\infty,H))/^tG). \end{equation}
It suffices to show that this morphism induces an isomorphism on rational homology, modulo conjugation in the domain. In formulas, \eqref{eq:comp} sends a chain $H_0\leq \hdots \leq H_k$ (lying between $H$ and $G$) together with coordinates $(t_0,\hdots,t_k)\in \Delta^k$ to the loop
\begin{equation} \label{eq:1} [\widetilde{\sigma}^G((\alpha_2)_e(\{H_i\},\{t_i\}))\circ l_{|H|}]\in \Omega (S^1\wedge \overline{S}(\infty,H)). \end{equation}
Here, we have again used that the diagonal (in the sense of the last part of Section \ref{sec:geomidea}) formed with respect to the intervals $[\frac{l-1}{|H|},\frac{l}{|H|}]$ corresponds to precomposition with $l_{|H|}$, see Section \ref{sec:construction}. Multiplication by $l_{|H|}$ is a rational equivalence, and 
\[ \widetilde{\sigma}^G:\left(\Sigma^{\infty}\overline{S}(\infty,H)\right)/^tG\to \Omega \left(\left(\Sigma^{\infty}S^1\wedge \overline{S}(\infty,H)\right)/^tG\right) \]
is a stable equivalence (as it agrees with the map $\widetilde{\lambda}_{\left(\Sigma^{\infty}\overline{S}(\infty,H)\right)/^tG}$, see the beginning of Section \ref{sec:proof}). So we can leave both out of the formula. What remains is given by applying the suspension spectrum functor $\Sigma^{\infty}$ to a wedge of space level maps
\begin{align*}  \beta_H:L(G)^{[H,G]}/\{\text{non-max. chains}\}& \to \overline{S}(\infty,H) \\
							(H_i,t_i) & \mapsto	[(\alpha_2)_e(H_i,t_i)]=[p(\lr(\sum_{i=0}^k t_i\cdot e_{H_i}))]
\end{align*}
followed by the projection to the translation $G$-orbits 
\[ \bigvee \overline{S}(\infty,H)\to \left(\bigvee \overline{S}(\infty,H)\right)/^tG. \]
We first show:
\begin{Lemma} \label{lem:split} Each \[ \beta_H: L(G)^{[H,G]}/\{\textnormal{non-max. chains}\}\to \overline{S}(\infty,H) \]
is a weak equivalence of spaces. 
\end{Lemma}
\begin{proof} If $H$ is equal to $G$ and hence $n=1$, both sides consist of two points and $\beta_H$ is a bijection. So from now on we assume that $H$ is a proper subgroup of $G$ and consider the map
\begin{align} \label{eq:lemma} \overline{\beta}_H:L(G)^{[H,G]}& \to (S^{\infty \rrho_G})^H \\
              \nonumber (H_i,t_i)& \mapsto \lr(\sum_{i=0}^k t_i\cdot e_{H_i}).
\end{align}
The map $\beta_H$ in the statement of the lemma is obtained from $\overline{\beta}_H$ by quotienting out by all non-maximal chains in the domain and the subspace $A$ of the target given by all vectors which are either fixed by a larger subgroup than $H$ or have norm at least $\frac{1}{\sqrt{2}}$. Since both $L(G)^{[H,G]}$ and $(S^{\infty \rrho_G})^H$ are contractible and the inclusions of the respective subspaces are cofibrations, it suffices to show that $\overline{\beta}_H$ induces a weak equivalence between the complex of non-maximal chains and $A$. We note that the former is given by the pushout \[ L(G)^{(H,G]}\cup_{L(G)^{(H,G)}} L(G)^{[H,G)}\]
of half-closed respectively open subintervals of the subgroup lattice. The space $A$ can be expressed in a similar way: Let $A_1\leq A$ be the subspace of vectors that have norm at most $1$ and are fixed by a subgroup properly containing $H$, and $A_2\leq A$ the subspace of all vectors of length at least $\frac{1}{\sqrt{2}}$. Together, the two cover $A$.
Then $\overline{\beta}_H$ maps $L(G)^{(H,G]}$ into $A_1$ and, by the same proof as for Lemma \ref{lem:reduced}, $L(G)^{[H,G)}$ into $A_2$. Since all of the spaces $L(G)^{(H,G]},L(G)^{[H,G)},A_1$ and $A_2$ are contractible (the latter two can be contracted onto $0$ and $\infty$, respectively), we are left to show that $\overline{\beta}_H$ induces a weak equivalence from $L(G)^{(H,G)}$ to the intersection of $A_1$ and $A_2$. This intersection is the space of vectors of $(S^{\infty \rrho_G})^H$ that have norm in the interval $[\frac{1}{\sqrt{2}},1]$ and are fixed by a larger subgroup. It deformation retracts onto 
\[ S(G)^{(H,G)}\defeq \colim_{H\lneq K\lneq G} \left(S(\infty\cdot \overline{\rho}_G)^K\right). \] 
We now consider the following commutative square:
\[ \xymatrix{ L(G)^{(H,G)}=\colim_{H\lneq K\lneq G} \left(L(G)^{[K,G)}\right) \ar[r] & \colim_{H\lneq K\lneq G} \left(S(\infty \cdot \overline{\rho}_G)^K\right)=S(G)^{(H,G)} \\
	      \hocolim_{H\lneq K\lneq G} \left(L(G)^{[K,G)}\right) \ar[r] \ar[u] & \hocolim_{H\lneq K\lneq G} \left(S(\infty \cdot \overline{\rho}_G)^K\right) \ar[u]  }
\]
Both vertical maps are the respective canonical map from the homotopy colimit to the colimit, and the lower horizontal map is induced from the restriction of $\overline{\beta}_H$ to the intervals $[K,G)$. Since each $S(\infty \overline{\rho}_G)^K$ is contractible and so is $L(G)^{[K,G)}$, it follows by homotopy-invariance of homotopy colimits that the lower horizontal map is a weak equivalence. We claim that both vertical maps are also weak equivalences. For the left one, this can be seen by noting that it is split by the map 
\[ L(G)^{(H,G)}=\hocolim_{H\lneq K\lneq G} *\to \hocolim_{H\lneq K\lneq G} \left(L(G)^{[K,G)}\right) \]
induced from the inclusions $*\mapsto K\in (L(G)^{[K,G)})_0$. Again by homotopy-invariance of homotopy colimits, this splitting is a weak equivalence and hence so is the left vertical map. Finally, to derive that the right vertical map is a weak equivalence one can use a $G$-CW structure on $S(\infty \cdot \overline{\rho}_G)$ to apply a general statement: For every $G$-cell complex $X$ the colimit over the $X^K$ with $K\in (H,G)$ is also a homotopy colimit. This can be seen by checking it for the orbits $G/L$ and using that cell attachments are both colimits and homotopy colimits.
\end{proof}
Finally, we consider the following commutative diagram:
\[ \xymatrix{ H_*(\bigvee L(G)^{[H,G]}/\{\text{non-max. chains}\},\Q) \ar[r]^-{\cong} \ar[d] & H_*(\bigvee \overline{S}(\infty,H),\Q) \ar[d] \\
							(H_*(\bigvee L(G)^{[H,G]}/\{\text{non-max. chains}\},\Q))_{^cG} \ar[r]^-{\cong}	\ar[dr]_-{(*)}	& (H_*(\bigvee \overline{S}(\infty,H),\Q))_{^cG} \ar[d] \\
						  & H_*((\bigvee \overline{S}(\infty,H))/^tG,\Q) }
\]
We want to show that the map $(*)$ is an isomorphism. The middle horizontal map is an isomorphism, since the wedge of the $\beta_H$ is equivariant for the conjugation actions. Furthermore, $\bigvee \overline{S}(\infty,H)$ is a cofibrant $G$-space under the translation action, so there is a natural isomorphism 
\[ H_*(\left(\bigvee \overline{S}(\infty,H)\right)/^tG,\Q)\cong H_*(\bigvee \overline{S}(\infty,H),\Q)_{^tG}. \]
Hence it suffices to see that conjugation and translation induce the same action on rational homology. This comes out of the proof of Lemma \ref{lem:split} above: Up to equivalence we can replace $\bigvee \overline{S}(\infty,H)$ by the suspension of the unreduced suspension of
\[ \bigvee_H \left(\hocolim_{H\lneq K\lneq G} (S(\infty \cdot\overline{\rho}_G)^K)\right). \]
Collapsing all $S(\infty \cdot \overline{\rho}_G)^K$ to a point yields a weak equivalence to 
\[ \bigvee_H \left(\hocolim_{H\lneq K\lneq G} *\right)= \bigvee_H L(G)^{(H,G)}. \]
This equivalence turns both the translation and the conjugation action on the $S(\infty \cdot \overline{\rho}_G)^K$ into the conjugation action on the wedge of subgroup intervals. So the two actions agree on homology, which finishes the proof of Proposition \ref{prop:splitinj}.
\subsection{Equivalence of the subquotients} \label{sec:ad}
We now show by other means that the subquotients in our two filtrations have the same rational homotopy type, i.e., that there is a rational equivalence  \begin{equation} \label{eq:eqgeom} \Phi^G(Sp^n/Sp^{n-1})\simeq_{\Q}\Sigma^{\infty} (|\widetilde{L}(G)_n/\widetilde{L}(G)_{n-1}|)/G.\end{equation} Since the latter is a finite complex, it has finite dimensional rational homology. So our map
\[ \walph_*:(H_*(|\widetilde{L}_f(G)_n/\widetilde{L}_f(G)_{n-1}|,\Q))_G\to H_*(\Omega \shift \Phi^G(Sp^n/Sp^{n-1}),\Q) \]
from the previous section must also be surjective, hence an isomorphism, proving Theorem \ref{theo:strong}.

To deduce the equivalence \eqref{eq:eqgeom} we combine work of Arone \cite{Aro15}, Arone-Dwyer \cite{AD01},  Arone-Brantner \cite{AB18} and Schwede \cite{Sch17}. In \cite[Proposition 1.11]{Sch17}, Schwede showed that there is a $G$-stable equivalence
\[ Sp^n_G/Sp^{n-1}_G\simeq_G \Sigma^{\infty} (B_G{\mathcal{F}_n})^{\dia}. \]
Here, $B_G{\mathcal{F}_n}$ is $E_G{\mathcal{F}_n}/\Sigma_n$, where $E_G{\mathcal{F}_n}$ is a universal $G$-space for the family $\mathcal{F}_n$ of subgroups of $\Sigma_n$ that do not act transitively on $\underline{n}=\{1,\hdots, n\}$. This means that $E_G{\mathcal{F}_n}$ is a cofibrant $(\Sigma_n\times G)$-space with the following two properties:
\begin{enumerate}[(i)]
	\item The $\Sigma_n$-isotropy of every point in $E_G{\mathcal{F}_n}$ lies in $\F_n$.
	\item The fixed points $(E_G{\mathcal{F}_n})^H$ are contractible for every subgroup $H\leq \Sigma_n\times G$ for which the intersection $H\cap (\Sigma_n\times \{1\})$ lies in~$\F_n$.
\end{enumerate}
Furthermore, the superscript $(-)^{\dia}$ denotes the unreduced suspension of a space. The non-equivariant version of this statement was previously shown by Lesh in \cite{Les00}. In \cite[Section 7]{AD01}, also for the case~$G=\{e\}$, Arone and Dwyer gave another description of this suspension spectrum, which we now mimic in the equivariant context. For this we denote by $\Pi_n$ the $\Sigma_n$-poset of non-trivial proper partitions of the set $\underline{n}$, and by $E_G\Sigma_n$ a universal $G$-space for $\Sigma_n$, i.e., a cofibrant $(\Sigma_n\times G)$-space satisfying conditions (i) and (ii) above with $\mathcal{F}_n$ replaced by the family consisting only of the trivial subgroup of $\Sigma_n$. Our aim is to show:
\begin{Prop} \label{prop:ad} The based $G$-spaces $(B_G\F_n)^{\dia}$ and $(E_G\Sigma_n)_+\wedge_{\Sigma_n} (|\Pi_n|^{\dia}\wedge S^n)$ are $G$-weakly equivalent after one suspension. Hence, there is a $G$-stable equivalence
\[ Sp^n_G/Sp^{n-1}_G\simeq_G \Sigma^{\infty} ((E_G\Sigma_n)_+\wedge_{\Sigma_n} (|\Pi_n|^{\dia}\wedge S^n)) \]
\end{Prop}
Here, $\Sigma_n$ acts on $S^n$ by permuting the coordinates. The $G$-action on 
\[ (E_G\Sigma_n)_+\wedge_{\Sigma_n} (|\Pi_n|^{\dia}\wedge S^n) \]
is only through $E_G\Sigma_n$.

\begin{Remark} \label{rem:globalquot} There also exists a global version of this result. Let $E_{gl}\Sigma_n$ denote a global universal space for $\Sigma_n$. For example, in the framework of orthogonal spaces of \cite[Section 1.1]{Sch18} or \cite[Section 1]{Sch17}, a model for $E_{gl}\Sigma_n$ is given by $L(\R^n,-)$ with $\Sigma_n$ permuting the coordinates of $\R^n$. Then there is a global equivalence
\[ Sp^n/Sp^{n-1}\simeq_{gl} \Sigma^{\infty} ((E_{gl}\Sigma_n)_+\wedge_{\Sigma_n} (|\Pi_n|^{\dia}\wedge S^n)). \]
The underlying $G$-space of $E_{gl}\Sigma_n$ is $E_G\Sigma_n$, so this gives back Proposition \ref{prop:ad} for all finite $G$. 
The proof we give below would also work in the global setting and prove this stronger statement, but we refrain from doing so to ease the exposition.
\end{Remark}

The following two pages are devoted to proving Proposition \ref{prop:ad}. The arguments are equivariant adaptions of the arguments in \cite[Section 7]{AD01}. Let $\Sing(E_G{\mathcal{F}_n})$ denote the subspace of $E_G{\mathcal{F}_n}$ consisting of the points with non-trivial $\Sigma_n$-isotropy. We note that $\Sing(E_G{\mathcal{F}_n})$ is a universal $G$-space for $\F_n^{\circ}$, the collection of non-transitive subgroups minus the trivial subgroup. Then there is a $(\Sigma_n\times G)$-cofiber sequence
\[ E_G{\Sigma_n}_+\wedge (\Sing(E_G{\mathcal{F}_n}))^{\dia}\to (\Sing(E_G{\mathcal{F}_n}))^{\dia} \to (E_G\Sigma_n*\Sing(E_G{\mathcal{F}_n}))^{\dia}\simeq_{\Sigma_n\times G}(E_G\mathcal{F}_n)^{\dia}, \]
where $*$ denotes the join. Smashing with the $\Sigma_n$-equivariant diagonal inclusion $i:S^1\to S^n$ yields a commutative diagram of $(\Sigma_n\times G)$-spaces
\[ \xymatrix{ E_G{\Sigma_n}_+\wedge \Sing(E_G{\mathcal{F}_n})^{\dia}\wedge S^1 \ar[r] \ar[d] & \Sing(E_G{\mathcal{F}_n})^{\dia}\wedge S^1 \ar[r] \ar[d] & (E_G\mathcal{F}_n)^{\dia}\wedge S^1 \ar[d] \\
E_G{\Sigma_n}_+\wedge\Sing(E_G{\mathcal{F}_n})^{\dia}\wedge S^n \ar[r] & \Sing(E_G{\mathcal{F}_n})^{\dia}\wedge S^n\ar[r]  &  (E_G\mathcal{F}_n)^{\dia}\wedge S^n.}
\]
We have:
\begin{Lemma} \label{lem:ad1} The map \[ (E_G\mathcal{F}_n)^{\dia}\wedge i:(E_G\mathcal{F}_n)^{\dia}\wedge S^1\to (E_G\mathcal{F}_n)^{\dia}\wedge S^n\] is a based $(\Sigma_n\times G)$-homotopy equivalence.
\end{Lemma}
\begin{Lemma} \label{lem:ad2} The quotient $(\Sing(E_G\F_n)^{\dia}\wedge S^n)/\Sigma_n$ is based $G$-weakly contractible. 
\end{Lemma}
The proofs of these lemmas are given below. Hence, quotienting out $\Sigma_n$ in the lower cofiber sequence yields a cofiber sequence of based $G$-spaces, which by Lemma \ref{lem:ad2} exhibits $((E_G\F_n)^{\dia}\wedge S^n)/\Sigma_n$ as the suspension of $((E_G\Sigma_n)_+\wedge_{\Sigma_n} (\Sing(E_G\F_n)^{\dia}\wedge S^n))$. By Lemma \ref{lem:ad1}, the former is $G$-weakly equivalent to $(B_G\F_n)^{\dia}\wedge S^1$, so we obtain:
\begin{Cor} There are $G$-weak equivalences
\begin{align*} (B_G\F_n)^{\dia}\wedge S^1 & \simeq ((E_G\Sigma_n)_+\wedge_{\Sigma_n} (\Sing(E_G\F_n)^{\dia}\wedge S^n))\wedge S^1 \\
& \simeq ((E_G\Sigma_n)_+\wedge_{\Sigma_n} ((E_G\F_n^{\circ})^{\dia}\wedge S^n))\wedge S^1. \end{align*}
\end{Cor}

\begin{proof}[Proof of Lemma \ref{lem:ad1}] Since both sides are cofibrant based $(\Sigma_n\times G)$-spaces, it suffices to show that the map induces a weak equivalence on all fixed point spaces. Let $H$ be a subgroup of $\Sigma_n\times G$. If the intersection $H\cap (\Sigma_n\times \{1\})$ acts non-transitively on $\underline{n}$, the $H$-fixed points of $(E_G\mathcal{F}_n)^{\dia}$ are contractible and hence the map is necessarily a weak equivalence. If $H\cap (\Sigma_n\times \{1\})$ does act transitively, the inclusion $i^H:S^1\to (S^n)^H$ is even a homeomorphism (in fact for this it would suffice that the projection of $H$ to $\Sigma_n$ act transitively). So $(E_G\mathcal{F}_n)^{\dia}\wedge i$ induces a weak equivalence on fixed points for all subgroups of $\Sigma_n\times G$, which finishes the proof.
\end{proof}

\begin{proof}[Proof of Lemma \ref{lem:ad2}] For this we make use of the specific model for $E_G\F_n$ that comes out of Schwede's proof \cite[Proposition 1.11]{Sch17}. It is given by $S(\overline{\R^n}\otimes \U_G)$, the unit sphere in the tensor product of the reduced natural $\Sigma_n$-representation with a complete $G$-universe. What we need from this model is the property that all the $\Sigma_n$-isotropy lies in \emph{complete subgroups}, i.e., subgroups of $\Sigma_n$ that are conjugate to one of the form $\Sigma_{n_1}\times \hdots \times \Sigma_{n_k}$ with all $n_i>0$, $\sum n_i=n$ and $k>1$. Hence, the isotropy of $\Sing(S(\overline{\R^n}\otimes \U_G))$ lies in non-trivial complete subgroups.

We now prove that more generally, the quotient $(X\wedge S^n)/\Sigma_n$ is $G$-contractible for any based cofibrant $(\Sigma_n\times G)$-space $X$ with all $\Sigma_n$-isotropy non-trivial and complete, or possibly the whole group $\Sigma_n$ (which we have to include because of the cone points of $\Sing(E_G\Pi_n)^{\dia}$). Without loss of generality we can assume that $X$ is a $(\Sigma_n\times G)$-cell complex. Via induction over the cells and passing to the sequential colimit we can reduce to showing that $((\Sigma_n\times G)/H_+\wedge A\wedge S^n)/\Sigma_n$ is $G$-weakly contractible for any space $A$ with trivial $(\Sigma_n\times G)$-action and $H\leq \Sigma_n\times G$ a subgroup for which $H\cap (\Sigma_n\times \{1\})$ is non-trivial and complete or equal to $\Sigma_n$.
We denote this intersection by $H'$ and the projection of $H$ to $G$ by $K$. For every $k\in K$ we choose an element $\psi(k)\in \Sigma_n$ such that $(\psi(k),k)$ lies in $H$. This property uniquely characterizes $\psi(k)$ up to multiplication with an element in $H'$, and every $\psi(k)$ automatically lies in the normalizer of $H'$. Altogether, $k\mapsto [\psi(k)]$ defines a homomorphism $\overline{\psi}:K\to W_{\Sigma_n}H'$ into the Weyl group. 
Then there is a $G$-homeomorphism
\[ ((\Sigma_n\times G)/H_+\wedge A\wedge S^n)/\Sigma_n \cong G\ltimes_K (A\wedge (S^n/H')), \]
with $K$ acting on $S^n/H'$ via restriction along $\overline{\psi}$.
So it suffices to see that $S^n/H'$ is $(W_{\Sigma_n}H')$-equivariantly contractible. Up to conjugacy, $H'$ is of the form 
\[ \Sigma_{n_1}^{\times i_1}\times \hdots \times \Sigma_{n_k}^{\times i_k} \]
with all $n_j$ pairwise different and $\sum (i_j \cdot n_j)=n$. The Weyl group is given by $\Sigma_{i_1}\times \hdots \times \Sigma_{i_k}$. Then, $S^n/H'$ is homeomorphic to \[ (S^{n_1}/\Sigma_{i_1})^{\wedge i_1}\wedge \hdots \wedge (S^{n_k}/\Sigma_{n_k})^{\wedge i_k},\] with the Weyl group permuting the smash factors in each $(S^{n_j}/\Sigma_{n_j})^{\wedge i_j}$. By \cite[Lemma 7.10]{AD01}, $(S^{n_j}/\Sigma_{n_j})$ is contractible whenever $n_j$ is greater than $1$. Since $H'$ is non-trivial, this has to be the case for some $n_j$, which finishes the proof.
\end{proof}

\begin{Remark} A more conceptual way to phrase the first part of the proof of Lemma \ref{lem:ad2} would be to say that the universal $G$-space for the collection of complete subgroups is $(\Sigma_n\times G)$-weakly equivalent to the universal $G$-space for the family of non-transitive subgroups $E_G\F_n$. This follows directly from the fact that $S(\overline{\R^n}\otimes \U_G)$ is a universal $G$-space for both, but could also be proved along the lines of \cite[Lemma 4.3]{AD01}.
\end{Remark}

One can further simplify $E_G\F_n^{\circ}$. For this we think of $\F_n^{\circ}$ as a $\Sigma_n$-poset, ordered by inclusion.

\begin{Lemma} There is a $(\Sigma_n\times G)$-map $E_G\F_n^{\circ}\to |\F_n^{\circ}|$ that induces a $(\Sigma_n\times G)$-weak equivalence
\[ E_G\Sigma_n\times E_G\F_n^{\circ} \xr{\simeq} E_G\Sigma_n \times |\F_n^{\circ}|. \]
\end{Lemma}
\begin{proof} This is a general fact about universal $G$-spaces for collections. Another model for $E_G\F_n^{\circ}$ (see \cite[Section 2]{AD01}, for example) is given by the nerve of the following category $\mathcal{E}_G\F_n^{\circ}$: An object $(M,x)$ is a transitive $(\Sigma_n\times G)$-set $M$ with $\Sigma_n$-isotropy in $\F_n^{\circ}$ together with a chosen element $x\in M$, and morphisms are given by $(\Sigma_n\times G)$-equivariant maps that preserve the chosen elements. The $(\Sigma_n\times G)$-action fixes the $(\Sigma_m\times G)$-sets and permutes the element $x$. This category carries a $(\Sigma_n\times G)$-functor $F$ to $\F_n^{\circ}$ by sending $(M,x)$ to the $\Sigma_n$-isotropy of $x$.

We claim that this functor does the job. For this we have to check that the fixed point functor \[ F^{\Gamma(\alpha)}:(\mathcal{E}_G\F_n^{\circ})^{\Gamma(\alpha)}\to (\F_n^{\circ})^{\Gamma(\alpha)}=(\F_n^{\circ})^{\im(\alpha)}\] induces a weak equivalence on nerves for every graph of a group homomorphism $\alpha$ from a subgroup $H$ of $G$ to $\Sigma_n$. But this follows from the fact that each such $F^{\Gamma(\alpha)}$ has a left adjoint, given by sending an $\im(\alpha)$-fixed subgroup $K$ of $\Sigma_n$ to $(\Sigma_n\times G/\langle K,\Gamma(\alpha)\rangle, [1])$, where $\langle K,\Gamma(\alpha)\rangle$ denotes the subgroup generated by $K\times \{1\}$ and $\Gamma(\alpha)$.
\end{proof}

Finally we relate $\F_n^{\circ}$ to the partition poset $\Pi_n$. There is a $\Sigma_n$-map of posets $j:\Pi_n\to \F_n^{\circ}$ sending a partition $\underline{n}=M_1\sqcup \hdots \sqcup M_k$ to the associated non-transitive subgroup $\Sigma_{M_1}\times \hdots \times \Sigma_{M_k}$.
\begin{Lemma} The map $j:\Pi_n\to \F_n^{\circ}$ induces a $\Sigma_n$-weak equivalence after applying the nerve.
\end{Lemma}
\begin{proof} An equivariant left adjoint is given by the map of posets that associates to every non-transitive and non-trivial subgroup $H\leq \Sigma_n$ the partition of $\underline{n}$ into $H$-orbits.
\end{proof}

This finishes the proof of Proposition \ref{prop:ad}.

Since there are natural isomorphisms $\Phi^G (\Sigma^{\infty} X)\cong \Sigma^{\infty}X^G$ for based $G$-spaces $X$, Proposition \ref{prop:ad} reduces the computation of $\Phi^G(Sp^n/Sp^{n-1})$ to the computation of the $G$-fixed points of the $G$-space $(E_G\Sigma_n)_+\wedge _{\Sigma_n} (|\Pi_n|^{\dia}\wedge S^n)$. To determine these, we make use of the following lemma.

\begin{Lemma} Let $K$ and $G$ be finite groups, and $X$ a cofibrant based $(K\times G)$-space such that the $K$-action is free away from the basepoint. Then there is a homeomorphism 
\[ (X/K)^G\cong \bigvee_{(\alpha:G\to K)} X^{\Gamma(\alpha)}/C(\alpha), \]
where the wedge is taken over conjugacy classes of group homomorphisms $\alpha:G\to K$, $C(\alpha)\leq K$ is the centralizer of the image of $\alpha$ and $\Gamma(\alpha)\leq K\times G$ is the graph of $\alpha$.
\end{Lemma}
\begin{proof} A proof for the version where $X$ is unbased (and not necessarily cofibrant) can be found in {\cite[Proposition B.17]{Sch18}}. The same arguments show that in the based case the canonical continuous map from the wedge of the $X^{\Gamma(\alpha)}/C(\alpha)$ to $(X/K)^G$ is bijective. The cofibrancy assumption ensures that it is a homeomorphism.
\end{proof}

Since $G$ acts trivially on $|\Pi_n|^{\dia}\wedge S^n$, application of this lemma gives a homeomorphism
\[ ((E_G\Sigma_n)_+\wedge _{\Sigma_n} (|\Pi_n|^{\dia}\wedge S^n))^G \cong \bigvee_{(\alpha:G\to \Sigma_n)} EC(\alpha)_+\wedge_{C(\alpha)} ((|\Pi_n|^{\dia}\wedge S^n)^{\im(\alpha)})
\]
Each $\alpha:G\to \Sigma_n$ defines a $G$-set structure on $\underline{n}$. The centralizer $C(\alpha)$ is given by the automorphisms of that $G$-set, and the fixed points $(S^n)^{\im(\alpha)}$ are a sphere of dimension the number of its $G$-orbits. So, written in a more coordinate free way, we obtain a homeomorphism
\[ ((E_G\Sigma_n)_+\wedge _{\Sigma_n} (|\Pi_n|^{\dia}\wedge S^n))^G\cong \bigvee_{(M\ G\text{-set},|M|=n)} (E{\Aut}_G(M))_+\wedge_{\Aut_G(M)} ((|\Pi_M|^G)^{\dia}\wedge (S^{M/G})),\]
where the wedge is taken over isomorphism classes of $G$-sets of order $n$. The fixed point sets $|\Pi_M|^G$ of partition posets that appear here have been studied by Arone \cite{Aro15} and Arone-Brantner \cite{AB18}. We can make use of their results to see that a lot of the wedge summands are (rationally) contractible, simplifying the expression for $((E_G\Sigma_n)_+\wedge _{\Sigma_n} (|\Pi_n|^{\dia}\wedge S^n))^G$. A finite $G$-set is called \emph{isotypical} if the isotropy groups of all of its points are conjugate.
\begin{Prop}[Fixed points of partition posets] We have:
\begin{enumerate}[(i)]
 \item If $M$ is not isotypical, the fixed points $|\Pi_M|^G$ are contractible and hence so is $(E{\Aut}_G(M))_+\wedge_{\Aut_G(M)} ((|\Pi_M|^G)^{\dia}\wedge (S^{M/G}))$.
 \item If $M$ is isotypical but not transitive, the space 
\[ (E{\Aut}_G(M))_+\wedge_{\Aut_G(M)} ((|\Pi_M|^G)^{\dia}\wedge (S^{M/G}))\]
is rationally contractible.
 \item If $M=G/H$ for a subgroup $H$ of $G$, then $\Pi_M^G$ is isomorphic to $L(G)^{(H,G)}$.
\end{enumerate}
\end{Prop}
\begin{proof}
$(i)$: This is \cite[Lemma 7.1]{Aro15}.

$(ii)$: Let $M=\bigsqcup_m G/H$ with $m\geq 2$. Then the automorphism group of $M$ is given by the wreath product $\Sigma_m\wr W_G(H)$. By Arone \cite[Proposition 9.1]{Aro15}, there is a $W_G(H)^m$-equivariant map
\begin{equation} \label{eq:ar} W_G(H)^m\ltimes_{W_G(H)}(|\Pi_m|*|\Pi_{G/H}^G|^{\dia})\xr{\simeq} |\Pi_M|^G \end{equation}
that is a non-equivariant equivalence, where $W_G(H)$ sits inside $W_G(H)^m$ diagonally. Direct inspection of its definition shows that the adjoint $(|\Pi_m|*|\Pi_{G/H}^G|^{\dia})\to |\Pi_M|^G$ (basically given by the cartesian product of partitions) is not only equivariant over the Weyl group, but also over the symmetric group $\Sigma_m$, if one lets it act on $\Pi_M$ by permuting the $m$ copies of $G/H$. In other words, we can also think of $\eqref{eq:ar}$ as a $(\Sigma_m\wr W_G(H))$-equivariant map
\[ (\Sigma_m\wr W_G(H))\ltimes_{\Sigma_m\times W_G(H)}(|\Pi_m|*|\Pi_{G/H}^G|^{\dia})\xr{\simeq} |\Pi_M|^G \]
that is a non-equivariant equivalence. Hence, we find that 
\[ E\Aut(M)_+\wedge_{\Aut(M)} ((|\Pi_M|^G)^{\dia}\wedge S^{M/G}) \]
is weakly equivalent to
\[ ((E\Sigma_m)_+\wedge_{\Sigma_m} (|\Pi_m|^{\dia}\wedge S^m))\wedge (EW_G(H)_+\wedge_{W_G(H)} (|\Pi_{G/H}|^G)^{\dia})\wedge S^1. \]
Here, we used that $|\Pi_m|*|\Pi_{G/H}^G|^{\dia}\simeq |\Pi_m|^{\dia}\wedge |\Pi_{G/H}^G|^{\dia}$, as described in \cite[Lemma 2.5]{Aro15}. So it suffices to note that $(E\Sigma_m)_+\wedge_{\Sigma_m} (|\Pi_m|^{\dia}\wedge S^m)$ is rationally contractible for $m\geq2$. This can for example be seen by replacing $|\Pi_m|$ by the weakly equivalent $E\Pi_m$ (the universal space for the collection of complete subgroups of $\Sigma_n$) and using that the strict quotient $(E\Pi_m^{\dia}\wedge S^m)/\Sigma_m$ is contractible. The latter was shown in the proof of Lemma \ref{lem:ad2}.

Part $(iii)$ follows from the fact that a $G$-fixed partition of $G/H$ is determined by its summand containing $H/H$ and that this summand has to be of the form $K/H$ for a subgroup $H\lneq K\lneq G$ (see \cite[Lemma 7.2]{Aro15}).
\end{proof}

The automorphism group of $G/H$ is given by the Weyl group $W_G(H)$, and so we see that there is a rational equivalence
\begin{align*} ((E_G\Sigma_n)_+\wedge_{\Sigma_n}(|\Pi_n|^{\dia}\wedge S^n))^G & \simeq_{\Q} \bigvee_{\substack{(H\leq G) \\ [G:H]=n }} (EW_G(H))_+\wedge_{W_G(H)}(|L(G)^{(H,G)}|^{\dia}\wedge S^1)
\end{align*}
We claim that the right hand side is rationally equivalent to $|\widetilde{L}_n(G)/\widetilde{L}_{n-1}(G)|/G$. In fact, we already saw in the proof of Lemma \ref{lem:split} that $|\widetilde{L}_n(G)/\widetilde{L}_{n-1}(G)|$ is isomorphic to the wedge over all index $n$ subgroups $H$ of the spaces \begin{equation*} \label{eq:quot} |L(G)^{[H,G]}|/(|L(G)^{[H,G)}|\cup_{|L(G)^{(H,G)}|} |L(G)^{(H,G]}|)\simeq |L(G)^{(H,G)}|^{\dia}\wedge S^1.\end{equation*}
After taking $G$-orbits, we can equivalently form the wedge over representatives of conjugacy classes of such subgroups $H$, and quotient each summand by the Weyl-group action. Since orbits and homotopy orbits are rationally equivalent for finite groups, the claim follows.

So together with Proposition \ref{prop:ad} and the fact that geometric fixed points commute with suspension spectra, this shows that there is a rational equivalence 
\[\Phi^G(Sp^n/Sp^{n-1})\simeq_{\Q}\Sigma^{\infty} (|\widetilde{L}(G)_n/\widetilde{L}(G)_{n-1}|/G), \]
which finishes the proof of our main result, Theorem \ref{theo:strong}.

\section{Examples} \label{sec:examples}
In this section we go through some small examples, where we describe the respective filtrations on subgroup lattices and read off the resulting rational equivariant homotopy groups of the symmetric products. In all cases we discuss, the quotient $L(G)/G$ is isomorphic to the nerve of the poset of conjugacy classes of subgroups of $G$. In general there is a natural surjection from the former to the latter, which is not always injective, as one can see - for example - with the symmetric group $\Sigma_5$: The transpositions $(12)$ and $(34)$ are conjugate in~$\Sigma_5$, but not in the subgroup $\Sigma_2\times \Sigma_3$. Since $\Sigma_2\times \Sigma_3$ is its own normalizer, this implies that the chains $\langle (12)\rangle \leq \Sigma_2\times \Sigma_3$ and $\langle (34) \rangle \leq \Sigma_2\times \Sigma_3$ represent different $1$-simplices in $L(\Sigma_5)/\Sigma_5$, while they represent the same $1$-simplex in the nerve of the poset of conjugacy classes.

\begin{Example}[Symmetric group $\Sigma_3$] We start with the symmetric group on $3$ letters. On the left we depict the subgroup lattice modulo conjugation, on the right the filtration by the $L(\Sigma_3)_n/\Sigma_3$ and the resulting dimensions for the rational $\Sigma_3$-homotopy groups of the symmetric products. All rational homotopy groups~$\pi_k^{\Sigma_3}(Sp^n)\otimes \Q$ with $k>1$ are trivial.

\newcommand{\unit}{0.4}
\renewcommand{\arraystretch}{1.3}

\center{\begin{tabular} {c c c | c}
\xymatrix@=2mm{
					& \Sigma_3 & \\
					& & A_3\ar@{-}[lu]_{2} \\
					\Sigma_2 \ar@{-}[ruu]^{3} & & \\
					& \{e\} \ar@{-}[lu]^2 \ar@{-}[uur]_3 &}
&&& \\
\end{tabular}
\begin{tabular} { >{$}c<{$} | >{$}c<{$} | >{$}c<{$} | >{$}c<{$} | >{$}c<{$}}
n & 1 & 2 & [3,5] & [6,\infty] \\ \hline
L(\Sigma_3)_n/\Sigma_3 &

\raisebox{-.43\height}{\begin{tikzpicture}
\coordinate (e) at (0,0);
\coordinate (s2) at (-\unit,\unit);
\coordinate (a3) at (\unit,2*\unit);
\coordinate (s3) at (0,3*\unit);

\foreach \x in {e,s2,a3,s3}{
\filldraw[black] (\x) circle (1pt); }

\filldraw[white] (0,3.4*\unit) (0.1pt);
\filldraw[white] (0,-.4*\unit) (0.1pt);
\end{tikzpicture}}

& 

\raisebox{-.42\height}{\begin{tikzpicture}
\coordinate (e) at (0,0);
\coordinate (s2) at (-\unit,\unit);
\coordinate (a3) at (\unit,2*\unit);
\coordinate (s3) at (0,3*\unit);

\draw[black] (e) -- (s2);
\draw[black] (a3) -- (s3);

\foreach \x in {e,s2,a3,s3}{
\filldraw[black] (\x) circle (1pt); }

\end{tikzpicture}}

&

\raisebox{-.42\height}{\begin{tikzpicture}
\coordinate (e) at (0,0);
\coordinate (s2) at (-\unit,\unit);
\coordinate (a3) at (\unit,2*\unit);
\coordinate (s3) at (0,3*\unit);

\draw[black] (e) -- (s2) -- (s3) -- (a3) -- cycle;

\foreach \x in {e,s2,a3,s3}{
\filldraw[black] (\x) circle (1pt); }

\end{tikzpicture}}

& 

\raisebox{-.42\height}{\begin{tikzpicture}
\coordinate (e) at (0,0);
\coordinate (s2) at (-\unit,\unit);
\coordinate (a3) at (\unit,2*\unit);
\coordinate (s3) at (0,3*\unit);

\draw[black,fill=yellow!50] (e) -- (s2) -- (s3) -- (a3) -- cycle;
\draw[black] (e) -- (s2) -- (s3) -- (a3) -- cycle;

\foreach \x in {e,s2,a3,s3}{
\filldraw[black] (\x) circle (1pt); }

\end{tikzpicture}}

\\ \hline
\dim(\pi_1^{\Sigma_3}(Sp^n)\otimes \Q) & 0 & 0 & 1 & 0 \\ 
\dim(\pi_0^{\Sigma_3}(Sp^n)\otimes \Q) & 4 & 2 & 1 & 1 
\end{tabular} }
\end{Example}

\begin{Example}[Dihedral group $D_8$] For the dihedral group with $16$ elements the filtration stabilizes, up to homotopy, at $n=4$. All minimal subgroup inclusions are of index $2$.

\newcommand{\unit}{0.5}
\renewcommand{\arraystretch}{1.3}
\center{\begin{tabular}{c c c | c}
\xymatrix@=1.5mm{& D_8 & \\
							D_4^{(1)} \ar@{-}[ur] & C_8 \ar@{-}[u] & D_4^{(2)} \ar@{-}[ul] \\
							D_2^{(1)} \ar@{-}[u] & C_4 \ar@{-}[u] \ar@{-}[ur] \ar@{-}[ul] & D_2^{(2)} \ar@{-}[u] \\
							D_1^{(1)} \ar@{-}[u] & C_2 \ar@{-}[u] \ar@{-}[ur] \ar@{-}[ul] & D_1^{(2)} \ar@{-}[u] \\
							& \{e\} \ar@{-}[u] \ar@{-}[ur] \ar@{-}[ul] &}
& &&\\
\end{tabular}
\begin{tabular} {c | c | c | c}
$n$ & $1$ & $2,3$ & $[4,\infty]$ \\ \hline
$L(D_8)_n/D_8$

&

\raisebox{-0.42\height}{\begin{tikzpicture}
\coordinate (e) at (0,0);
\coordinate (c2) at (0,\unit);
\coordinate (c4) at (0,2*\unit);
\coordinate (c8) at (0,3*\unit);
\coordinate (d8) at (0,4*\unit);
\coordinate (d11) at (-\unit,\unit);
\coordinate (d12) at (\unit,\unit);
\coordinate (d21) at (-\unit,2*\unit);
\coordinate (d22) at (\unit,2*\unit);
\coordinate (d31) at (-\unit,3*\unit);
\coordinate (d32) at (\unit,3*\unit);

\foreach \x in {e,c2,c4,c8,d8,d11,d12,d21,d22,d31,d32}{
\filldraw[black] (\x) circle (1pt); };
\filldraw[white] (0,4.3*\unit,0) circle (0.1pt);
\filldraw[white] (0,-.25*\unit,0) (0.1pt);
\end{tikzpicture}}

&
\raisebox{-0.42\height}{\begin{tikzpicture}
\coordinate (e) at (0,0);
\coordinate (c2) at (0,\unit);
\coordinate (c4) at (0,2*\unit);
\coordinate (c8) at (0,3*\unit);
\coordinate (d8) at (0,4*\unit);
\coordinate (d11) at (-\unit,\unit);
\coordinate (d12) at (\unit,\unit);
\coordinate (d21) at (-\unit,2*\unit);
\coordinate (d22) at (\unit,2*\unit);
\coordinate (d41) at (-\unit,3*\unit);
\coordinate (d42) at (\unit,3*\unit);

\draw[black] (e) -- (c2) -- (c4) -- (c8) -- (d8);
\draw[black] (e) -- (d11) -- (d21) -- (d41) -- (d8);
\draw[black] (e) -- (d12) -- (d22) -- (d42) -- (d8);
\draw[black] (c2) -- (d21);
\draw[black] (c2) -- (d22);
\draw[black] (c4) -- (d41);
\draw[black] (c4) -- (d42);

\foreach \x in {e,c2,c4,c8,d8,d11,d12,d21,d22,d41,d42}{
\filldraw[black] (\x) circle (1pt); };
\end{tikzpicture}}

&

\raisebox{-0.42\height}{\begin{tikzpicture}
\coordinate (e) at (0,0);
\coordinate (c2) at (0,\unit);
\coordinate (c4) at (0,2*\unit);
\coordinate (c8) at (0,3*\unit);
\coordinate (d8) at (0,4*\unit);
\coordinate (d11) at (-\unit,\unit);
\coordinate (d12) at (\unit,\unit);
\coordinate (d21) at (-\unit,2*\unit);
\coordinate (d22) at (\unit,2*\unit);
\coordinate (d41) at (-\unit,3*\unit);
\coordinate (d42) at (\unit,3*\unit);

\draw[black,fill=yellow!50] (e) -- (d11) -- (d41) -- (d8) -- cycle;
\draw[black,fill=yellow!50] (e) -- (d12) -- (d42) -- (d8) -- cycle;

\draw[black] (c2) -- (d21);
\draw[black] (c2) -- (d22);
\draw[black] (c4) -- (d41);
\draw[black] (c4) -- (d42);

\foreach \x in {e,c2,c4,c8,d8,d11,d12,d21,d22,d41,d42}{
\filldraw[black] (\x) circle (1pt); };
\end{tikzpicture}}

\\ \hline

$\dim(\pi_1^{D_8}(Sp^n)\otimes \Q)$ & $0$ & $6$ & $0$ \\

$\dim(\pi_0^{D_8}(Sp^n)\otimes \Q)$ & $11$ & $1$ & $1$ \\
\end{tabular} }
\end{Example}
\begin{Example}[Special linear group $SL_2(\mathbb{F}_3)$] The figure below depicts the filtration for the special linear group $SL_2(\mathbb{F}_3)$, the semi-direct product of the quaternion group $Q_8$ with the cyclic group $C_3$. As in the previous examples, the bottom
two rows in the table give the dimensions of the rationalized $\pi_0^{SL_2(\mathbb{F}_3)}(Sp^n)$ and $\pi_1^{SL_2(\mathbb{F}_3)}(Sp^n)$.

\newcommand{\unit}{0.4}
\renewcommand{\arraystretch}{1.3}

\center{
\begin{tabular}{c | c}
\xymatrix@=1.5mm{& SL_2(\mathbb{F}_3) &\\
								C_6 \ar@{-}[ur]^-{4} & &  Q_8 \ar@{-}[ul]_-{3} \\
								C_3 \ar@{-}[u]^{2} & & C_4 \ar@{-}[u]_{2}\\
								  & & C_2 \ar@{-}[uull]_3 \ar@{-}[u]_{2} \\
								& \{e\} \ar@{-}[uul]^{3} \ar@{-}[ur]_{2} & &} &
\end{tabular}								
\begin{tabular} {c | c | c | c | c | c}
$1$ & $2$ & $3$ & $4,5$ & $[6,11]$ & $[12,\infty]$ \\ \hline 				
\raisebox{-0.422\height}{\begin{tikzpicture}
\coordinate (t) at (0,4.5*\unit);
\coordinate (1) at (0,0);
\coordinate (c3) at (-\unit,2*\unit);
\coordinate (c2) at (\unit,\unit);
\coordinate (c4) at (\unit,2*\unit);
\coordinate (c6) at (-\unit,3*\unit);
\coordinate (q8) at (\unit,3*\unit);
\coordinate (s) at (0,4*\unit);

\foreach \x in {1,c3,c2,c4,c6,q8,s}{
\filldraw[black] (\x) circle (1pt); }

\filldraw[white] (t) circle (0.1pt);
\filldraw[white] (0,-.4*\unit) circle (0.1pt);

\end{tikzpicture}}
&
\raisebox{-0.42\height}{\begin{tikzpicture}

\coordinate (1) at (0,0);
\coordinate (c3) at (-\unit,2*\unit);
\coordinate (c2) at (\unit,\unit);
\coordinate (c4) at (\unit,2*\unit);
\coordinate (c6) at (-\unit,3*\unit);
\coordinate (q8) at (\unit,3*\unit);
\coordinate (s) at (0,4*\unit);

\draw[black] (1) -- (c2);
\draw[black] (c2) -- (c4);
\draw[black] (c4) -- (q8);
\draw[black] (c3) -- (c6);

\foreach \x in {1,c3,c2,c4,c6,q8,s}{
\filldraw[black] (\x) circle (1pt); }
\end{tikzpicture}}
&

\raisebox{-0.42\height}{\begin{tikzpicture}

\coordinate (1) at (0,0);
\coordinate (c3) at (-\unit,2*\unit);
\coordinate (c2) at (\unit,\unit);
\coordinate (c4) at (\unit,2*\unit);
\coordinate (c6) at (-\unit,3*\unit);
\coordinate (q8) at (\unit,3*\unit);
\coordinate (s) at (0,4*\unit);

\draw[black] (1) -- (c2) -- (c4) -- (q8) -- (s);
\draw[black] (1) -- (c3) -- (c6) -- (c2);

\foreach \x in {1,c3,c2,c4,c6,q8,s}{
\filldraw[black] (\x) circle (1pt); }
\end{tikzpicture}}
&
\raisebox{-0.42\height}{\begin{tikzpicture}

\coordinate (1) at (0,0);
\coordinate (c3) at (-\unit,2*\unit);
\coordinate (c2) at (\unit,\unit);
\coordinate (c4) at (\unit,2*\unit);
\coordinate (c6) at (-\unit,3*\unit);
\coordinate (q8) at (\unit,3*\unit);
\coordinate (s) at (0,4*\unit);

\draw[black] (1) -- (c2) -- (c4) -- (q8) -- (s) -- (c6);
\draw[black] (1) -- (c3) -- (c6) -- (c2);

\foreach \x in {1,c3,c2,c4,c6,q8,s}{
\filldraw[black] (\x) circle (1pt); }

\end{tikzpicture}}
&
\raisebox{-0.42\height}{\begin{tikzpicture}

\coordinate (1) at (0,0);
\coordinate (c3) at (-\unit,2*\unit);
\coordinate (c2) at (\unit,\unit);
\coordinate (c4) at (\unit,2*\unit);
\coordinate (c6) at (-\unit,3*\unit);
\coordinate (q8) at (\unit,3*\unit);
\coordinate (s) at (0,4*\unit);

\draw[black] (c2) -- (c4) -- (q8) -- (s) -- (c6);
\draw[black,fill=yellow!50] (1) -- (c3) -- (c6) -- (c2)-- cycle;
\draw[black] (1) -- (c3) -- (c6) -- (c2)-- cycle;

\foreach \x in {1,c3,c2,c4,c6,q8,s}{
\filldraw[black] (\x) circle (1pt); }

\end{tikzpicture}}
&
\raisebox{-0.42\height}{\begin{tikzpicture}

\coordinate (1) at (0,0);
\coordinate (c3) at (-\unit,2*\unit);
\coordinate (c2) at (\unit,\unit);
\coordinate (c4) at (\unit,2*\unit);
\coordinate (c6) at (-\unit,3*\unit);
\coordinate (q8) at (\unit,3*\unit);
\coordinate (s) at (0,4*\unit);

\draw[black,fill=green!70] (c2) -- (c4) -- (q8) -- (s) -- (c6) -- cycle;
\draw[black,fill=yellow!50] (1) -- (c3) -- (c6) -- (c2)-- cycle;
\draw[black] (c2) -- (c4) -- (q8) -- (s) -- (c6) -- cycle;
\draw[black] (1) -- (c3) -- (c6) -- (c2)-- cycle;

\foreach \x in {1,c3,c2,c4,c6,q8,s}{
\filldraw[black] (\x) circle (1pt); }

\end{tikzpicture}} \\ \hline
$0$ & $0$ & $1$ & $2$ & $1$ & $0$ \\ 
$7$ & $3$ & $1$ & $1$ & $1$ & $1$\\
\end{tabular}
}

\end{Example}
\begin{Example}[Cyclic groups] For cyclic groups $C_m$, the subgroups correspond to divisors of $m$ ordered by divisibility. If $m$ is the product of $k$ different primes, the resulting lattice is a $k$-dimensional cube. In particular, the subcomplex $L(C_m)_{m-1}$ is the boundary of the cube and hence isomorphic to $S^{k-1}$. This shows that for every $k\in \N$ there exists an $n\in \N$ and a finite group $G$ such that $\pi_k^G(Sp^n)\otimes \Q$ is non-trivial. Below is the $3$-dimensional cube for the example $m=30=2\cdot 3\cdot 5$, where the elements at the top are those divisible by $2$, the ones on the right those divisible by $3$ and the ones at the back those divisible by $5$. The bottom three rows in the table list the resulting dimensions of the rationalized $\pi^{C_{30}}_0(Sp^n)$, $\pi^{C_{30}}_1(Sp^n)$ and $\pi^{C_{30}}_2(Sp^n)$.

\newcommand{\Depth}{0.8}
\newcommand{\Height}{0.8}
\newcommand{\Width}{0.8}
\newcommand{\Diff}{2.5}
\renewcommand{\arraystretch}{1.3}

\center{\begin{tabular}{ c | c | c | c | c | c | c | c}
	$1$ & $2$ & $3,4$ & $5$ & $[6,9]$ & $[10,14]$ & $[15,29]$ & $[30,\infty]$ \\ \hline

\raisebox{-0.42\height}{\begin{tikzpicture}
\coordinate (t) at (0,1.2*\Width,0);
\coordinate (1) at (0,0,0);
\coordinate (2) at (0,\Width,0);
\coordinate (10) at (0,\Width,\Height);
\coordinate (5) at (0,0,\Height);
\coordinate (3) at (\Depth,0,0);
\coordinate (6) at (\Depth,\Width,0);
\coordinate (30) at (\Depth,\Width,\Height);
\coordinate (15) at (\Depth,0,\Height);

\foreach \x in {1, 2, 3, 5, 6, 10, 15, 30}{
\filldraw[black] (\x) circle (1pt); } 
\filldraw[white] (t) circle (0.1pt);
\filldraw[white] (0,-0.55*\Width,0) circle (0.1pt);
\end{tikzpicture}}

&

\raisebox{-0.42\height}{\begin{tikzpicture}

\coordinate (1a) at (0+\Diff,0,0);
\coordinate (2a) at (0+\Diff,\Width,0);
\coordinate (10a) at (0+\Diff,\Width,\Height);
\coordinate (5a) at (0+\Diff,0,\Height);
\coordinate (3a) at (\Depth+\Diff,0,0);
\coordinate (6a) at (\Depth+\Diff,\Width,0);
\coordinate (30a) at (\Depth+\Diff,\Width,\Height);
\coordinate (15a) at (\Depth+\Diff,0,\Height);

\foreach \x in {1, 2, 3, 5, 6, 10, 15, 30}{
\filldraw[black] (\x a) circle (1pt); } 

\draw[black] (1a) -- (2a);
\draw[black] (3a) -- (6a);
\draw[black] (5a) -- (10a);
\draw[black] (15a) -- (30a);

\end{tikzpicture}}

&

\raisebox{-0.42\height}{\begin{tikzpicture}
\coordinate (1b) at (0+2*\Diff,0,0);
\coordinate (2b) at (0+2*\Diff,\Width,0);
\coordinate (10b) at (0+2*\Diff,\Width,\Height);
\coordinate (5b) at (0+2*\Diff,0,\Height);
\coordinate (3b) at (\Depth+2*\Diff,0,0);
\coordinate (6b) at (\Depth+2*\Diff,\Width,0);
\coordinate (30b) at (\Depth+2*\Diff,\Width,\Height);
\coordinate (15b) at (\Depth+2*\Diff,0,\Height);

\foreach \x in {1, 2, 3, 5, 6, 10, 15, 30}{
\filldraw[black] (\x b) circle (1pt); } 

\draw[black] (1b) -- (2b);
\draw[black] (3b) -- (6b);
\draw[black] (5b) -- (10b);
\draw[black] (15b) -- (30b);
\draw[black] (1b) -- (3b);
\draw[black] (2b) -- (6b);
\draw[black] (5b) -- (15b);
\draw[black] (10b) -- (30b);

\end{tikzpicture}}

&

\raisebox{-0.42\height}{\begin{tikzpicture}
\coordinate (1c) at (0+3*\Diff,0,0);
\coordinate (2c) at (0+3*\Diff,\Width,0);
\coordinate (10c) at (0+3*\Diff,\Width,\Height);
\coordinate (5c) at (0+3*\Diff,0,\Height);
\coordinate (3c) at (\Depth+3*\Diff,0,0);
\coordinate (6c) at (\Depth+3*\Diff,\Width,0);
\coordinate (30c) at (\Depth+3*\Diff,\Width,\Height);
\coordinate (15c) at (\Depth+3*\Diff,0,\Height);

\foreach \x in {1, 2, 3, 5, 6, 10, 15, 30}{
\filldraw[black] (\x c) circle (1pt); } 

\draw[black] (1c) -- (2c);
\draw[black] (3c) -- (6c);
\draw[black] (5c) -- (10c);
\draw[black] (15c) -- (30c);
\draw[black] (1c) -- (3c);
\draw[black] (2c) -- (6c);
\draw[black] (5c) -- (15c);
\draw[black] (10c) -- (30c);
\draw[black] (1c) -- (5c);
\draw[black] (2c) -- (10c);
\draw[black] (3c) -- (15c);
\draw[black] (6c) -- (30c);

\end{tikzpicture}}

&

\raisebox{-0.42\height}{\begin{tikzpicture}

\coordinate (1d) at (0+4*\Diff,0,0);
\coordinate (2d) at (0+4*\Diff,\Width,0);
\coordinate (10d) at (0+4*\Diff,\Width,\Height);
\coordinate (5d) at (0+4*\Diff,0,\Height);
\coordinate (3d) at (\Depth+4*\Diff,0,0);
\coordinate (6d) at (\Depth+4*\Diff,\Width,0);
\coordinate (30d) at (\Depth+4*\Diff,\Width,\Height);
\coordinate (15d) at (\Depth+4*\Diff,0,\Height);

\foreach \x in {1, 2, 3, 5, 6, 10, 15, 30}{
\filldraw[black] (\x d) circle (1pt); } 

\draw[black] (1d) -- (5d);
\draw[black,fill=yellow!50] (1d) -- (2d) -- (6d) -- (3d) -- cycle;% Back Face
\draw[black,fill=yellow!50,opacity=0.8] (5d) -- (10d) -- (30d) -- (15d) -- cycle;% Front Face
\draw[black] (3d) -- (15d);
\draw[black] (2d) -- (10d);
\draw[black] (6d) -- (30d);

\filldraw[black] (30d) circle (1pt);

\end{tikzpicture}}

&

\raisebox{-0.42\height}{\begin{tikzpicture}

\coordinate (1e) at (0+5*\Diff,0,0);
\coordinate (2e) at (0+5*\Diff,\Width,0);
\coordinate (10e) at (0+5*\Diff,\Width,\Height);
\coordinate (5e) at (0+5*\Diff,0,\Height);
\coordinate (3e) at (\Depth+5*\Diff,0,0);
\coordinate (6e) at (\Depth+5*\Diff,\Width,0);
\coordinate (30e) at (\Depth+5*\Diff,\Width,\Height);
\coordinate (15e) at (\Depth+5*\Diff,0,\Height);

\foreach \x in {1, 2, 3, 5, 6, 10, 15, 30}{
\filldraw[black] (\x e) circle (1pt); } 

\draw[black,fill=yellow!50] (1e) -- (2e) -- (6e) -- (3e) -- cycle;% Back Face
\draw[black,fill=green!50] (1e) -- (2e) -- (10e) -- (5e) -- cycle;% Left Face
\draw[black,fill=yellow!50,opacity=0.8] (5e) -- (10e) -- (30e) -- (15e) -- cycle;% Front Face
\draw[black,fill=green!70,opacity=0.8] (3e) -- (6e) -- (30e) -- (15e) -- cycle;% Right Face

\filldraw[black] (30e) circle (1pt);

\end{tikzpicture}}

&

\raisebox{-0.42\height}{\begin{tikzpicture}

\coordinate (1f) at (0+6*\Diff,0,0);
\coordinate (2f) at (0+6*\Diff,\Width,0);
\coordinate (10f) at (0+6*\Diff,\Width,\Height);
\coordinate (5f) at (0+6*\Diff,0,\Height);
\coordinate (3f) at (\Depth+6*\Diff,0,0);
\coordinate (6f) at (\Depth+6*\Diff,\Width,0);
\coordinate (30f) at (\Depth+6*\Diff,\Width,\Height);
\coordinate (15f) at (\Depth+6*\Diff,0,\Height);

\foreach \x in {1, 2, 3, 5, 6, 10, 15, 30}{
\filldraw[black] (\x f) circle (1pt); } 

\draw[black,fill=red!50] (1f) -- (3f) -- (15f) -- (5f) -- cycle;% Bottom Face
\draw[black,fill=yellow!50] (1f) -- (2f) -- (6f) -- (3f) -- cycle;% Back Face
\draw[black,fill=green!50] (1f) -- (2f) -- (10f) -- (5f) -- cycle;% Left Face
\draw[black,fill=yellow!50,opacity=0.8] (5f) -- (10f) -- (30f) -- (15f) -- cycle;% Front Face
\draw[black,fill=green!70,opacity=0.8] (3f) -- (15f) -- (30f) -- (6f) -- cycle;% Right Face
\draw[black,fill=red!50,opacity=0.8] (2f) -- (6f) -- (30f) -- (10f) -- cycle;% Top Face

\filldraw[black] (30f) circle (1pt);

\end{tikzpicture}}

&

\raisebox{-0.42\height}{\begin{tikzpicture}

\coordinate (1g) at (0+7*\Diff,0,0);
\coordinate (2g) at (0+7*\Diff,\Width,0);
\coordinate (10g) at (0+7*\Diff,\Width,\Height);
\coordinate (5g) at (0+7*\Diff,0,\Height);
\coordinate (3g) at (\Depth+7*\Diff,0,0);
\coordinate (6g) at (\Depth+7*\Diff,\Width,0);
\coordinate (30g) at (\Depth+7*\Diff,\Width,\Height);
\coordinate (15g) at (\Depth+7*\Diff,0,\Height);

\foreach \x in {1, 2, 3, 5, 6, 10, 15, 30}{
\filldraw[black] (\x g) circle (1pt); }

\draw[black,fill=yellow!50] (5g) -- (10g) -- (30g) -- (15g) -- cycle;
\draw[black,fill=green!70] (3g) -- (15g) -- (30g) -- (6g) -- cycle;
\draw[black,fill=red!50] (2g) -- (6g) -- (30g) -- (10g) -- cycle;

\filldraw[black] (30g) circle (1pt);

\end{tikzpicture}}

\\ \hline $0$ & $0$ & $0$ & $0$ & $0$ & $0$ & $1$ & $0$\\
$0$ & $0$ & $2$ & $5$ & $3$ & $1$ & $0$ & $0$\\
$8$ & $4$ & $2$ & $1$ & $1$ & $1$ & $1$ & $1$\\
\end{tabular} }

\end{Example}
We close this section with a few general remarks on the complexes of the form $L(G)_{|G|-1}$, the last non-trivial stage in the filtration. They can be identified with the unreduced suspension of the lattice $L(G)^{(1,G)}$ of proper non-trivial subgroups of $G$, about which there are various results in the literature. For example, it was shown by Kratzer-Th{\'e}venaz \cite{KT85} that if $G$ is solvable, then $L(G)^{(1,G)}$ is homotopy-equivalent to a wedge of spheres of dimension two less than the chief length $c(G)$ of $G$. By Th{\'e}venaz \cite{The85}, the top homology $H_{c(G)-2}(L(G)^{(1,G)},\Z)$ is a permutation representation under the conjugation $G$-action.
The co-invariants $(H_{c(G)-2}(L(G)^{(1,G)},\Z))_G$ are still acted on by the outer automorphism group $\Out(G)$ of $G$. As we saw (Theorem \ref{thm:mackey}), after tensoring with $\Q$ this corresponds to the action of $\Out(G)$ on $\pi_{c(G)-1}^G(Sp^{|G|-1})\otimes \Q$ which is part of the structure of a global functor. These actions can be interesting representation-theoretically:
\begin{Example}[$(\Z/p)^n$ and the Steinberg module] \label{exa:steinberg} When $G=(\Z/p)^n$, the complex $L(G)^{(1,G)}$ is the Tits building for $\Out(G)=GL_n(\mathbb{F}_p)$. So, by a theorem of Solomon (\cite{Sol69}),  its homology $H_{n-2}(L(G)^{(1,G)},\Q)$ (and hence also $\pi_{n-1}^{G}(Sp^{p^n-1})\otimes \Q$) is isomorphic to the rational Steinberg module, a distinguished irreducible $GL_n(\mathbb{F}_p)$-representation of dimension $p^{\frac{n(n-1)}{2}}$. For example, when $p=n=2$, the Steinberg module is the reduced natural representation of $\Sigma_3\cong GL_2(\mathbb{F}_2)$.

A different relation between symmetric products of spheres and the Steinberg module - over $\mathbb{F}_p$ instead of $\Q$ - played a major role in Arone-Dwyer \cite{AD01}.
\end{Example}

\section{Global properties of $Sp^n_{\Q}$}
\label{sec:globalprop}
In this final section we describe homological properties of the $\Out^{op}$-complex models $\cC \widetilde{L}_n$ for the rational symmetric products. We first show that they are degreewise projective and then use this to prove that for $1<n<\infty$ they are not formal, i.e., not quasi-isomorphic to their homology with trivial differential. As a consequence, the rationalization $Sp^n_{\Q}$ is not a product of Eilenberg-MacLane spectra for any $n$ except $1$ and $\infty$. This is a truly global phenomenon, since over a fixed finite group $G$ every rational $G$-spectrum is determined by its homotopy groups. Finally, we give a proof that the cyclic $p$-groups are the only groups for which $\pi_*^G(Sp^n)\otimes \Q$ is always concentrated in degree $0$.

\begin{Prop} Each $\cC \widetilde{L}_n$ is degreewise projective as a $\Q[\Out^{op}]$-module.
\end{Prop}
\begin{proof} We need the following notion: A chain of subgroup inclusions $H_0\leq \hdots \leq H_k$ is called \emph{simple} if $H_0$ does not contain a non-trivial normal subgroup of $H_k$, i.e., if this chain cannot be obtained via pull-back along a surjective group homomorphism that is not an isomorphism.

Recall that the $k$-th level of $\mathcal{C} \widetilde{L}_n(G)$ is given by the $\Q$-linearization of the set of conjugacy classes of chains of proper subgroup inclusions which end in $G$, are of length $k$ and have total index at most $n$. The map 
\begin{align*}
\left(\bigsqcup_{H_0\lneq \hdots \lneq H_k \text{ simple}, [H_k:H_0]\leq n} \Epi(G,H_k)\right)/\text{iso} & \to  \{\text{$k$-chains of index $\leq n$ in $G$}\} \\
(H_0,\lneq \hdots, \lneq H_k,\psi:G\twoheadrightarrow H_k) & \mapsto  (\psi^{-1}(H_0)\lneq \hdots \lneq \psi^{-1}(H_k)=G)
\end{align*}
defines a natural bijection. On the left hand side, two pairs 
\[ (H_0\lneq \hdots \lneq H_k;\psi:G\twoheadrightarrow H_k) \]
and 
\[ (H'_0\lneq \hdots \lneq H'_k;\psi':G\twoheadrightarrow H'_k) \]
are considered isomorphic if there exists an isomorphism $\varphi:H_k\xr{\cong} H_k'$ that takes the first chain to the second and which satisfies $\varphi\circ \psi=\psi'$. For the inverse map one associates to a $k$-chain $H_0\lneq \hdots \lneq H_k=G$ the simple $k$-chain 
\[ H_0/H\lneq \hdots \lneq H_k/H=G/H \]
together with the projection $G\twoheadrightarrow G/H$, where $H$ is the intersection of all $G$-conjugates of $H_0$, or in other words the largest subgroup of $H_0$ that is normal in $G$. Hence, modding out by conjugations on both sides and restricting the sum to representatives of isomorphism classes, we find that there is an isomorphism of $\Out^{op}$-modules
\begin{equation} \label{eq:decomp} (\mathcal{C}  \widetilde{L}_n)_k\cong \bigoplus_{[H_0\lneq \hdots \lneq H_k \text{ simple}, [H_k:H_0]\leq n]} \Q[\Out(-,H_k)]/\Out(H_0\lneq \hdots \lneq H_k). \end{equation}
Here, $\Out(H_0\lneq \hdots \lneq H_k)$ denotes the group of conjugacy classes of automorphisms of $H_k$ which map the chain $H_0\lneq \hdots \lneq H_k$ to a conjugate of itself. The modules $\Q[\Out(-,H_k)]$ are by definition representable, hence projective. Furthermore, the orbits of a projective module under any action of a finite group $K$ are again projective, since the projection is split by the map $[x]\mapsto \frac{1}{|K|}\sum_{k\in K} (k\cdot x)$. This finishes the proof.
\end{proof}

The proof also applies to $n=\infty$ and hence $\cC \widetilde{L}$ itself, showing that it gives a projective resolution of the $\Out^{op}$-module that sends the trivial group to $\Q$ and all other finite groups to $0$. Using once more that the functor $\tau$ of Section \ref{sec:geom} is an equivalence, this shows that $\cC L$ is a projective resolution of the constant global functor $\underline{\Q}$.

We can use our algebraic model to see:

\begin{Prop} \label{prop:formal} For $1<n<\infty$ the rationalization $Sp^n_{\Q}$ is not a product of global Eilenberg-MacLane spectra. 
\end{Prop}
\begin{proof} Using Theorem \ref{theo:strong}, the statement follows if we show that the $\Out^{op}$-complex $\mathcal{C}\wL_n$ is not quasi-isomorphic to its homology with trivial differential.

Each $\mC \wL_n$ is concentrated in finitely many degrees $0,\hdots,\lfloor \log_2(n) \rfloor=a(n)$. We show that the highest possible $k$-invariant is non-trivial, i.e., that the map 
\[ \Sigma^{a_n}H_{a(n)}(\mC \wL_n)\to \mC \wL_n\]
does not have a section in the derived category. Since $\mC \wL(n)$ is degreewise projective, this is equivalent to the inclusion $H_{a(n)}(\mC \wL_n)\hookrightarrow (\mC \wL_n)_{a(n)}$ of $\Out^{op}$-modules not having a section. In fact we claim that any $\Out^{op}$-map $(\mC \wL_n)_{a(n)}\to H_{a(n)}(\mC \wL_n)$ is necessarily zero on all abelian groups. To see this we use the decomposition \eqref{eq:decomp} of $(\mC \wL_n)_{a(n)}$ above. When restricted to abelian groups, the only summands that play a role are those associated to chains $H_0\lneq \hdots \lneq H_{a(n)}$ with $H_{a(n)}$ abelian, since there is no surjective map from an abelian group to a non-abelian one. In the abelian case, the simpleness implies that $H_0$ is the trivial group, and so the order of $H_{a(n)}$ is 
at 
most $n$.
So, over abelian groups, $\mC \wL(n)$ is a quotient of a direct sum of representables for groups $G$ of order at most $n$. It now suffices to see that any map from these representables to $H_{a(n)}(\mC \wL_n)$ is trivial. By the Yoneda Lemma, such maps correspond to homology classes $H_{a(n)}(\mC \wL_n)(G)$. Since $\widetilde{L}(G)_n$ is contractible for groups $G$ of order $\leq n$ (unless $G$ is the trivial group, which can only appear if $a(n)=0$ and hence $n=1$), this homology is trivial. This proves the claim.

Hence it suffices to show that there exists an abelian group $G$ for which $H_{a(n)}(\mC \wL_n)(G)$ is non-trivial. By Example \ref{exa:steinberg}, such a $G$ is given by $(\Z/2)^{a(n)+1}$.
\end{proof}

Using similar arguments for other classes of groups, one can show the non-vanishing of many more $k$-invariants of $\mC \wL_n$. 

\begin{Remark} In contrast, the rationalized subquotients $(Sp^n/Sp^{n-1})_{\mathbb{Q}}$ are \emph{always} products of global Eilenberg-MacLane spectra. One can see this as follows: As noted in Remark \ref{rem:globalquot}, the global homotopy type of $Sp^n/Sp^{n-1}$ is given by 
\[ \Sigma^{\infty}((E_{gl}\Sigma_n)_+\wedge_{\Sigma_n} (|\Pi_n|^{\dia}\wedge S^n)), \]
or written differently 
\[ (E_{gl}\Sigma_n)_+\wedge_{\Sigma_n}(\Sigma^{\infty}(|\Pi_n|^{\dia}\wedge S^n)).\]
The construction $(E_{gl}\Sigma_n)_+\wedge_{\Sigma_n} (-)$ is a functor from $\Sigma_n$-spectra to global spectra that is left adjoint to the functor that takes a global spectrum to its underlying $\Sigma_n$-spectrum (see \cite[Theorem 4.5.24]{Sch18}). By considering its effect on geometric fixed points, it is not hard to see that this functor preserves rational Eilenberg-MacLane spectra. Furthermore, it commutes with rationalization. Hence, since any $\Sigma_n$-spectrum~$X$ decomposes into Eilenberg-MacLane spectra rationally (see Remark \ref{rem:gmorita}), it follows that so does any global spectrum of the form $(E_{gl}\Sigma_n)_+\wedge_{\Sigma_n} X$, and in particular~$Sp^n/Sp^{n-1}$.
\end{Remark}

Finally, we have:
\begin{Prop} \label{prop:onlypgroups} For a finite group $G$ the following are equivalent:
\begin{enumerate}[(i)]
	\item For all $n\in \N$ the graded vector space $\pi_*^G(Sp^n)\otimes \Q$ is concentrated in degree $0$.
	\item For all $n\in \N$ the vector space $\pi_1^G(Sp^n)\otimes \Q$ is trivial.
	\item $G\cong C_{p^n}$ for some prime $p$ and $n\in \N$.
\end{enumerate}
\end{Prop}
\begin{proof} If $G\cong C_{p^n}$, the subgroup lattice of $G$ is linear, and it is not hard to see that all subcomplexes $L(G)_n$ are either discrete (for $n<p$) or contractible. So, by Theorem \ref{thm:mackey}, $\pi_*^G(Sp^n)\otimes \Q$ is concentrated in degree $0$.

Item $(i)$ clearly implies $(ii)$, so it remains to show that $(ii)$ implies $(iii)$. For this we fix a finite group $G$ that is not cyclic of prime power order and show that some $\pi_1^G(Sp^n)\otimes \Q$ is non-trivial. To see this, we choose two non-conjugate maximal subgroups $H$ and $H'$ of $G$, which is possible since~$G$ is not cyclic of prime power order. For example, the existence of $H$ and $H'$ follows from the fact that the union of all conjugates of a proper subgroup can never be all of $G$.
Let $n$ denote the index of $H$ in $G$, which we can without loss of generality assume to be at least as large as that of $H'$.
Then the formal difference 
\[ [H\leq G]-[H'\leq G] \in \mathcal{C} \wL_n(G)_1\]
is a non-trivial $1$-cycle. Since any proper subgroup of $H$ has index larger than $n$ in $G$ and $G$ is the only subgroup containing $H$, there are no non-degenerate $2$-simplices of $\wL (G)_n/G$ that have $H\leq G$ or any of its conjugates as a face. Hence, $[H\leq G]-[H'\leq G]$ cannot be a boundary and thus defines a non-trivial element in 
\[ H_1(\mathcal{C} \wL_n(G))\cong \Phi_1^G(Sp^n)\otimes \Q.\]
Since $\Phi_1^G(Sp^n)\otimes \Q$ is a quotient of $\pi_1^G(Sp^n)\otimes \Q$, the latter is also non-trivial, which finishes the proof.
\end{proof}

\bigskip
\footnotesize

  \textsc{Department of Mathematical Sciences, University of Copenhagen, Denmark}\par\nopagebreak
  \textit{E-mail address}: \texttt{hausmann@math.ku.dk}
\end{document}